\documentclass[11pt,letterpaper]{article}
\usepackage[T1]{fontenc}
\usepackage[utf8]{inputenc}
\usepackage[sc]{mathpazo}

\usepackage{amsmath,amssymb,amsthm,mathtools}
\usepackage[margin=1.05in,headheight=14pt,headsep=20pt,footskip=28pt]{geometry}
\usepackage{microtype}
\usepackage{booktabs,array,enumitem,needspace}
\usepackage{titlesec,fancyhdr,etoolbox}
\usepackage{tikz}
\usetikzlibrary{positioning,cd}
\usepackage{hyperref}
\usepackage[nameinlink,noabbrev]{cleveref}
\definecolor{citepink}{RGB}{214,51,132}
\hypersetup{colorlinks=true,linkcolor=black,urlcolor=black,citecolor=citepink,
 pdftitle={Constructions and Upper Bounds for Trifferent Codes},
 pdfsubject={Trifferent codes, constructions, binary boxes},
 bookmarksnumbered=true,pdfdisplaydoctitle=true}
\titleformat{\section}{\fontsize{14}{17}\selectfont\bfseries}{\thesection}{0.7em}{}
\titleformat{\subsection}{\normalsize\bfseries}{\thesubsection}{0.65em}{}
\titleformat{\subsubsection}{\normalsize\itshape}{\thesubsubsection}{0.65em}{}
\titlespacing*{\section}{0pt}{3ex plus .7ex minus .2ex}{1.2ex plus .2ex}
\titlespacing*{\subsection}{0pt}{2.4ex plus .5ex minus .2ex}{1ex plus .2ex}
\setlist{itemsep=2pt,topsep=4pt,leftmargin=1.6em}
\allowdisplaybreaks[1]
\theoremstyle{plain}
\newtheorem{theorem}{Theorem}[section]
\newtheorem{lemma}[theorem]{Lemma}
\newtheorem{proposition}[theorem]{Proposition}

\theoremstyle{definition}
\newtheorem{definition}[theorem]{Definition}

\newtheorem{remark}[theorem]{Remark}
\newcommand{\F}{\mathbb F}
\newcommand{\E}{\mathbb E}

\newcommand{\word}[1]{\texttt{#1}}

\newcommand{\gstar}{1.162579}
\newcommand{\paperauthors}{%
	\begin{tabular}{c@{\hspace{3em}}c}
		Siddharth Bhandari & Abhishek Khetan\\
		{\small\texttt{siddharthbhandari.engg@gmail.com}} & {\small\texttt{abhishek.khetan@ashoka.edu.in}}
	\end{tabular}%
}
\newcommand{\paperaffiliations}{}
\fancypagestyle{plain}{\fancyhf{}\fancyfoot[C]{\small\thepage}}
\begin{document}
\thispagestyle{plain}
\begin{center}
{\fontsize{21}{25}\selectfont\bfseries Constructions and Upper Bounds\par for Trifferent Codes}\par
\ifdefempty{\paperauthors}{}{\vspace{10pt}{\normalsize\paperauthors}\par}
\ifdefempty{\paperaffiliations}{}{\vspace{4pt}{\small\paperaffiliations}\par}
\end{center}
\begin{abstract}
A ternary code of blocklength $n$ is trifferent if every three distinct codewords have a coordinate at which they take all three values.
Let $T(n)$ be the largest size of such a code.
We prove that $T(n)=\Omega(1.162579^n)$ and $T(n)=O\bigl(n^{-1/2}(3/2)^n\bigr)$.
The lower bound improves the exponential base $(9/5)^{1/4}=1.15829\ldots$ of K\"orner and Marton.
It comes from a random system of overlapping local constraints, in which every position is read two or three times, followed by thinning and the deletion of nearby pairs and of the remaining triples that violate trifference.
The required estimates reduce to finitely many inequalities, which we verify rigorously by computer.
The upper bound improves our earlier bound $O\bigl(n^{-2/5}(3/2)^n\bigr)$, using a compatibility constraint between the empty boxes attached to overlapping six-codeword configurations.
\end{abstract}

\setcounter{tocdepth}{2}
\tableofcontents


\section{Introduction}
A \emph{ternary word of blocklength $n$} is an element of $\{\word{0},\word{1},\word{2}\}^n$.
A coordinate $i\in\{1,\ldots,n\}$ is \emph{rainbow} for three ternary words $x,y,z$ of blocklength $n$ if $\{x_i,y_i,z_i\}=\{\word{0},\word{1},\word{2}\}$.
A set $C$ of ternary words of blocklength $n$ is a \emph{trifferent code} if every three distinct words of $C$ have a rainbow coordinate.
Let $T(n)$ denote the largest size of a trifferent code of blocklength $n$.
The \emph{Hamming distance} $d(x,y)$ of two words $x,y$ of blocklength $n$ is the number of coordinates $i$ with $x_i\ne y_i$.
A code has \emph{minimum distance at least} $\Delta$ if every two distinct words of the code are at Hamming distance at least $\Delta$.
A sequence of codes $C_n$ of blocklength $n$ has \emph{asymptotic relative minimum distance at least} $\gamma$ if $C_n$ has minimum distance at least $(\gamma-o(1))n$ as $n\to\infty$.
Our results concern two different parts of this extremal problem.

\begin{theorem}\label{thm:three}
	\[
		T(n)
		=\Omega(1.162579^n),
		\qquad
		T(n)=O\bigl(n^{-1/2}(3/2)^n\bigr).
	\]
\end{theorem}

The two bounds are proved by independent arguments, which we now outline.

To construct a large trifferent code we use the alteration method.
Call a triple of words \emph{bad} if it has no rainbow coordinate.
We choose a random set of words and then delete some of them, so that no bad triple remains.
The resulting code is large if the random set is large and contains few bad triples.
Our random set is produced as follows.
We choose at random a system of overlapping local constraints on words over a nine-letter alphabet, which are encoded as ternary words through the tetracode \cite{KornerMarton1988}, and we keep each word that satisfies every constraint independently with a fixed probability.
Choosing the system means choosing both how the constraints are wired to the variables and a permutation of the alphabet at each place where a constraint reads a variable, so there are three random choices, and they have different jobs.
The random permutations make the expected number of accepted words easy to compute.
The random wiring makes the expected numbers of pairs and triples of accepted words depend only on how many variables of each degree are in each equality state.
The final thinning balances the number of kept words against the cost of the deletions.
The constraints make the coordinates of these words dependent, and they are chosen so that many words satisfy them while close pairs and bad triples are rare.
Before deleting bad triples, we delete one word from every pair of kept words at small Hamming distance.
This removes, at the cost of a single word, every bad triple that contains such a pair.
All expectations are taken over the three random choices together, so the number of words satisfying one particular system never has to be controlled.
Generating functions express the expected counts, and weights that record distances reduce the required estimates to finitely many inequalities, which we verify by computer in exact rational arithmetic, with rigorous enclosures of the logarithms.

K\"orner and Marton obtained the construction base $(9/5)^{1/4}=1.15829\ldots$ using a nine-word inner code \cite{KornerMarton1988}.
Han, Sun and Ge subsequently improved its polynomial factor, obtaining $\Omega(\sqrt n(9/5)^{n/4})$ \cite{HanSunGe2026}.
Our construction improves the exponential base.
It is an existence proof based on a finite local rule and exact counting certificates, not a claim of an efficient deterministic encoder.

Yuan and Zhu \cite{YuanZhu2026} recently improved the K\"orner--Marton lower bound by sampling disjoint length-$14$ blocks independently over a nine-symbol alphabet, subject to a parity restriction within each block.
Our construction, developed independently before their preprint appeared, also uses a local restriction to make bad triples less likely.
The main difference is that our blocks overlap: each position is read two or three times by the local constraints.
Choosing the mix of degree-two and degree-three coordinates lets us control both close pairs and bad triples.
This gives an exponential base of $1.162579$, compared with approximately $1.158344$ in Yuan and Zhu's bound, and an asymptotic relative minimum distance of at least $9/32$.

To bound the size of an arbitrary trifferent code we work with binary boxes.
For a word $q$, the \emph{binary box} avoiding $q$ is the set of words $x$ with $x_i\ne q_i$ for every coordinate $i$.
Three words lie in a common binary box if and only if they have no rainbow coordinate, so a binary box contains at most two words of a trifferent code.
Each word lies in $2^n$ of the $3^n$ binary boxes, and averaging the bound of two words per box over all boxes gives the classical estimate $T(n)\le2(3/2)^n$.
Our improvement comes from comparing the boxes that contain two codewords with the boxes that contain none.
Using a compatibility constraint between the empty boxes attached to overlapping configurations of six codewords (Lemma~\ref{gen:coherence}), we show that, for codes with at least four words, the number of boxes containing two codewords is at most $O(1/n)$ times the number of empty boxes, and hence $O(3^n/n)$.
An estimate of the probability that two independent uniformly random codewords agree on a given set of coordinates then turns this count into the bound $T(n)=O(n^{-1/2}(3/2)^n)$.

For the upper bound, we improve our previous estimate \cite{BhandariKhetan2025} for arbitrary ternary trifferent codes from $O(n^{-2/5}(3/2)^n)$ to $O(n^{-1/2}(3/2)^n)$.
The improvement comes from studying how small configurations around empty binary boxes fit together.
This saves a factor of $n^{1/10}$ without changing the exponential base.

\paragraph{Proof roadmap.}
Section~\ref{lower:blocks} proves the lower bound.
Section~\ref{lower:sampling} recalls the sampling and deletion argument and the tetracode encoding, which lets us work with words over a nine-letter alphabet.
Section~\ref{lower:alteration} thins a random set instead of sampling from a law, and deletes close pairs before bad triples.
Section~\ref{lower:views} introduces systems of overlapping constraints, and Section~\ref{lower:counting} identifies two statistics of the local rule that govern the expected counts.
Section~\ref{lower:certified} fixes the local rule and the wiring, expresses the expected numbers of accepted words, pairs and bad triples through these statistics, and states the required counting estimates.
Section~\ref{lower:distance} reduces these estimates to finitely many inequalities and completes the proof.
Appendix~\ref{cert:appendix} describes the certificates for these inequalities and how they are checked in exact rational arithmetic with rigorous logarithm enclosures.
Section~\ref{sec:general} proves the upper bound, using the geometry of binary boxes and this compatibility constraint (Lemma~\ref{gen:coherence}).

Throughout, $[n]=\{1,\ldots,n\}$.
Logarithms are natural unless a base is displayed.
When samples are discussed, a bad triple may contain repeated words.
All code-size assertions concern sets of distinct words.
The supplementary programs verify finite certificates.
The analytic arguments below supply the conclusions for every blocklength.

\section{Construction}
\label{lower:blocks}\label{lower:source}\label{sec:preliminaries}

\begin{theorem}\label{lower:main}
	Let $g=1.162579$.
	At blocklengths $n=124k$ there are trifferent codes with $(1-o(1))g^n$ words.
	Consequently $T(n)=\Omega(g^n)$.
\end{theorem}

As a by-product of the proof, the codes at blocklengths $n=124k$ also have minimum distance at least $9n/32$.

\subsection{Sampling, deletion, and nine symbols}
\label{lower:sampling}

A natural way to build a trifferent code is the alteration method.
We draw a list of candidate words at random from a probability distribution, and then delete some of the words so that no bad triple remains.
Recall that a triple $(x,y,z)$ of words of blocklength $n$ is \emph{bad} if it has no rainbow coordinate, that is, if $|\{x_i,y_i,z_i\}|\le2$ for every $i\in[n]$.
The resulting code is large if the list is long and contains few bad triples.
We first make this precise for independent samples.

We identify the alphabet $\{\word{0},\word{1},\word{2}\}$ with the field $\F_3$, so the ternary words of blocklength $n$ are the elements of $\F_3^n$.
A \emph{law} on $\F_3^n$ is a probability distribution on $\F_3^n$.
The \emph{error} of a law $\mu$ on $\F_3^n$ is the probability that three independent samples from $\mu$ form a bad triple:
\[
 E_n(\mu)=\Pr[(X,Y,Z)\text{ is bad}],\qquad X,Y,Z\text{ independent, each with law }\mu.
\]
The three samples need not be distinct.
A triple in which two words coincide has no rainbow coordinate, so it is bad and contributes to $E_n(\mu)$.
Errors multiply over independent blocks (see below), so it is convenient to normalize by the blocklength.
The \emph{per-coordinate error} of a law $\mu$ on $\F_3^n$ is
\[
 \varepsilon(\mu)=E_n(\mu)^{1/n}.
\]

\begin{lemma}[Alteration]\label{lower:sampling-alteration}
	Let $\mu$ be a law on $\F_3^n$, and let $M\ge1$ be an integer with
	\[
		M-\binom M3E_n(\mu)\ge3.
	\]
	Then there is a trifferent code of blocklength $n$ with at least $M-\binom M3E_n(\mu)$ words.
\end{lemma}
\begin{proof}
Let $X_1,\ldots,X_M$ be independent samples from $\mu$.
For each set $\{i,j,k\}$ of three distinct indices, the triple $(X_i,X_j,X_k)$ is bad with probability $E_n(\mu)$.
Hence the expected number of index sets carrying a bad triple is $\binom M3E_n(\mu)$, and some outcome has at most this many.
Fix such an outcome, and delete one index from each index set carrying a bad triple.
At least $M-\binom M3E_n(\mu)\ge3$ indices remain, and no three of them carry a bad triple.
The words at the remaining indices are pairwise distinct.
Indeed, if two remaining indices carried the same word, then together with any third remaining index they would carry a bad triple.
So the remaining words form a set $C$ of at least $M-\binom M3E_n(\mu)$ words.
Any three distinct words of $C$ sit at three remaining indices, so they do not form a bad triple, and $C$ is trifferent.
\end{proof}

Lemma~\ref{lower:sampling-alteration} is applied most simply to laws that are repeated over independent blocks.
For laws $\mu$ on $\F_3^a$ and $\mu'$ on $\F_3^b$, we write $\mu\otimes\mu'$ for the law on $\F_3^{a+b}$ of the concatenation $(X,Y)$ of independent samples $X$ from $\mu$ and $Y$ from $\mu'$, and $\mu^{\otimes k}$ for the product of $k$ copies of $\mu$.
Take three independent samples from $\mu\otimes\mu'$.
Their first blocks are three independent samples from $\mu$, their second blocks are three independent samples from $\mu'$, and the two triples of blocks are independent.
A coordinate of the concatenation is rainbow if and only if it is a rainbow coordinate of one of the two triples of blocks.
So the three samples form a bad triple if and only if both triples of blocks are bad, and
\begin{equation}\label{eq:prelim-product}
 E_{a+b}(\mu\otimes\mu')=E_a(\mu)E_b(\mu').
\end{equation}
In particular, $E_{dk}(\mu^{\otimes k})=E_d(\mu)^k$ and $\varepsilon(\mu^{\otimes k})=\varepsilon(\mu)$ for every law $\mu$ on $\F_3^d$.

\begin{lemma}[Alteration for product laws]\label{prelim:alteration}
Let $\mu$ be a law on $\F_3^d$ with per-coordinate error $\varepsilon(\mu)<1$.
As $k\to\infty$, there are trifferent codes of blocklength $n=dk$ and size at least
\begin{equation}\label{eq:prelim-alteration}
 \left(\frac{2\sqrt2}{3}+o(1)\right)\varepsilon(\mu)^{-n/2}.
\end{equation}
\end{lemma}
\begin{proof}
Write $\varepsilon=\varepsilon(\mu)$.
By~\eqref{eq:prelim-product}, the law $\mu^{\otimes k}$ on $\F_3^n$ has error $E_d(\mu)^k=\varepsilon^n$.
Take $M=\lfloor\sqrt2\,\varepsilon^{-n/2}\rfloor$, so that $M^2\varepsilon^n\le2$.
Then $\binom M3\varepsilon^n\le M^3\varepsilon^n/6\le M/3$.
Since $\varepsilon<1$, we have $M\to\infty$, so $2M/3\ge3$ for all large $k$.
Lemma~\ref{lower:sampling-alteration} now gives a trifferent code of blocklength $n$ with at least $2M/3=(2\sqrt2/3+o(1))\varepsilon^{-n/2}$ words.
\end{proof}

So a law with per-coordinate error $\varepsilon(\mu)$ gives codes with growth $\varepsilon(\mu)^{-1/2}$ per coordinate.
Our construction refines this argument in two ways.
First, it deletes one word from every close pair before it deletes bad triples (Section~\ref{lower:alteration}).
Second, its candidates are not independent samples from a law.
They are the words accepted by a random system of constraints, each kept independently with a fixed probability (Sections~\ref{lower:alteration} and~\ref{lower:views}).

The construction works with words over a nine-letter alphabet and encodes them as ternary words at the end.
Let $\F_9$ be the field with nine elements.
We identify it with $\F_3^2$ as an $\F_3$-vector space (Appendix~\ref{lower:local-counts} gives an explicit model), call its elements \emph{labels}, and encode a label $(a,b)$ by
\[
 \psi(a,b)=(a,b,a+b,a-b)\in\F_3^4.
\]
The image $K=\psi(\F_3^2)$ is the nine-word tetracode, the inner code of K\"orner and Marton \cite{KornerMarton1988}.
The four coordinates of $\psi$ are the linear forms $a$, $b$, $a+b$ and $a-b$ on $\F_3^2$.
A direct check shows that exactly one of these four forms vanishes at each nonzero vector of $\F_3^2$.
Two consequences follow.
First, distinct words of $K$ differ in exactly three coordinates, since exactly one form vanishes on their difference.
Second, let $u,v,w\in\F_3^2$ be distinct.
Each of the differences $u-v$, $v-w$ and $w-u$ is annihilated by exactly one form, so at least one of the four forms vanishes on none of them.
That form takes three distinct values at $u$, $v$ and $w$, so the corresponding coordinate is rainbow for $\psi(u),\psi(v),\psi(w)$.
If instead two of $u,v,w$ are equal, then no coordinate is rainbow for $\psi(u),\psi(v),\psi(w)$.

For a word $x\in\F_9^m$ write $\Psi(x)=(\psi(x_1),\ldots,\psi(x_m))\in\F_3^{4m}$.
By the first consequence, the Hamming distance satisfies $d(\Psi(x),\Psi(y))=3\,d(x,y)$, where $d(x,y)$ is the number of positions at which the labels of $x$ and $y$ differ.
By the second, $(\Psi(x),\Psi(y),\Psi(z))$ is bad if and only if, at every position $j\in[m]$, two of the labels $x_j,y_j,z_j$ are equal, that is, if and only if $(x,y,z)$ is itself a bad triple of words over $\F_9$.
We therefore work with words over $\F_9$ until the final encoding step.

\subsection{Deleting close pairs first}
\label{lower:alteration}

Lemma~\ref{lower:sampling-alteration} deletes one word from every bad triple.
Many bad triples are bad for a simple reason: two of their words are equal or agree at most coordinates.
A triple with two equal words has no rainbow coordinate.
If two words of a triple agree at all but a few coordinates, a rainbow coordinate can occur only among those few coordinates.
We therefore first delete one word from every pair of words at small distance, and only afterwards delete one word from every bad triple that remains.
The first deletion removes every triple that contains such a pair, so the second deletion only has to deal with bad triples whose three words are pairwise far apart.
The first deletion also gives the code a large minimum distance.

Deleting pairs is cheap, because among $M$ words there are about $M^2/2$ pairs but $M^3/6$ triples.
If each pair is close with probability $P$, and each triple is bad and pairwise far apart with probability $E$, then the expected numbers of close pairs and of such triples are about $M^2P/2$ and $M^3E/6$, and the two deletions remove at most this many words in expectation.
This is small compared with $M$ if $MP$ and $M^2E$ are small, so $P$ may exceed $E$ by a factor of order $M$.
Such a gap is unavoidable when the words are drawn from a small set.
If $M$ words are drawn independently and uniformly from a set of $N$ words, any two of them coincide with probability $1/N$.
Since a triple with two equal words is bad, the uniform law on the set has error at least $1/N$, and the bound $M-\binom M3E_n$ of Lemma~\ref{lower:sampling-alteration} is at most of order $N^{1/2}$ for every $M$.
When close pairs are deleted first, coinciding words only contribute to $P$, and $M$ can be much larger than $N^{1/2}$.

We now make this precise for words over $\F_9$.
Unlike the independent sampling above, the construction thins a random set, so its candidates are distinct words.
Recall from Section~\ref{lower:sampling} that a triple $(x,y,z)$ of words of $\F_9^m$ is bad if at every position $j\in[m]$ two of the labels $x_j,y_j,z_j$ are equal, and that this holds if and only if $(\Psi(x),\Psi(y),\Psi(z))$ is bad.
Recall that $d(x,y)$ is the number of positions at which the labels of $x$ and $y$ differ.
Fix $\delta>0$.
A pair $(x,y)$ is \emph{close} if $d(x,y)<\delta m$.
In particular, $(x,x)$ is close.
A bad triple is \emph{far} if none of its three pairs is close.

Let $S$ be a random subset of $\F_9^m$, that is, a subset drawn from a probability distribution on the power set of $\F_9^m$.
Lemma~\ref{lower:two-stage} below holds for every such $S$.
From Section~\ref{lower:views} on, $S$ will be the set of words accepted by a random system of overlapping constraints.
Write $N=|S|$, and let $N_2$ be the number of close ordered pairs in $S^2$ and $N_3$ the number of far bad ordered triples in $S^3$.
Assume $\bar N=\E N>0$, and put
\begin{equation}\label{lower:annealed-counts}
 b_2=\frac{\E N_2}{\bar N^2},\qquad b_3=\frac{\E N_3}{\bar N^3}.
\end{equation}
These count ordered pairs and triples.
So $N_2$ includes the diagonal pairs $(x,x)$, and each unordered close pair or far bad triple is counted two or six times.
This costs only constant factors, and it makes $N_2$ and $N_3$ sums over assignments of states to the variables (Section~\ref{lower:distance}).

\begin{lemma}[Two-stage alteration]\label{lower:two-stage}
	For every real number $M$ with $0<M\le\bar N$, there is a set $C\subseteq\F_9^m$ with
	\begin{equation}\label{lower:joint-alteration}
		|C|\ge M-M^2b_2-M^3b_3
	\end{equation}
	such that $\Psi(C)$ is a trifferent code of blocklength $4m$ with minimum distance at least $3\delta m$.
\end{lemma}
\begin{proof}
	Put $p=M/\bar N$, so that $0<p\le1$.
	Given $S$, form the \emph{thinning} $T\subseteq S$ by keeping each word of $S$ independently with probability $p$.
	Let $P$ be the number of two-element subsets of $T$ whose two words form a close pair, and let $R$ be the number of three-element subsets of $T$ whose three words form a far bad triple.
	We take expectations over $S$ and $T$ together.

	First, $\E|T|=p\,\E N=M$.
	Next, condition on $S$.
	A two-element subset of $S$ survives in $T$ with probability $p^2$.
	The close two-element subsets of $S$ are at most $N_2$ in number, since each of them gives two close ordered pairs.
	So $\E P\le p^2\E N_2=M^2b_2$.
	In the same way, a far bad triple consists of three distinct words, since $d(x,x)=0<\delta m$.
	Each far bad three-element subset of $S$ survives with probability $p^3$, and these subsets are at most $N_3$ in number, so $\E R\le p^3\E N_3=M^3b_3$.
	Hence
	\[
		\E\bigl[|T|-P-R\bigr]\ge M-M^2b_2-M^3b_3.
	\]
	Fix an outcome of $S$ and $T$ in which $|T|-P-R$ is at least its expectation.

	Delete one word from every close two-element subset of $T$.
	Then delete one word from every far bad three-element subset of $T$ that is still intact.
	Let $C$ be the set of remaining words.
	At most $P+R$ words were deleted, so $|C|\ge|T|-P-R$, which gives~\eqref{lower:joint-alteration}.

	No two distinct words of $C$ form a close pair, so $d(x,y)\ge\delta m$ for all distinct $x,y\in C$.
	Since $d(\Psi(x),\Psi(y))=3\,d(x,y)$ (Section~\ref{lower:sampling}), the words of $\Psi(C)$ are distinct and have pairwise Hamming distance at least $3\delta m$.
	Now let $x,y,z\in C$ be distinct.
	Their three pairs are not close.
	If $(x,y,z)$ were bad, it would be far, and one of its words would have been deleted.
	So $(\Psi(x),\Psi(y),\Psi(z))$ is not bad, and $\Psi(C)$ is trifferent.
\end{proof}

The bound~\eqref{lower:joint-alteration} has the shape of Lemma~\ref{lower:sampling-alteration}.
The ratio $b_3$ plays the role of the error $E_n(\mu)$, and $b_2$ plays the role of the probability that two samples are close.
Indeed, if $N$ were always equal to $\bar N$, then $b_2$ and $b_3$ would be the probabilities that two or three independent uniform samples from $S$ form a close pair or a far bad triple.
The lemma is useful when $Mb_2$ and $M^2b_3$ are small.
Only far bad triples are counted, so the condition on triples is weaker than in Lemma~\ref{lower:sampling-alteration}.
The price is the condition on pairs, which involves only the first power of $M$.

The proof keeps each word of $S$ with a fixed probability $p$.
It does not sample words from the uniform law on $S$, as Lemma~\ref{lower:sampling-alteration} would.
Sampling from the uniform law would require, for a typical outcome of $S$ with $N>0$, the probability that three independent uniform samples from $S$ form a far bad triple.
This probability is $N_3/N^3$, a ratio of two random quantities, and it is not defined when $S$ is empty, which may happen with positive probability.
We can compute the expectations of counts such as $N$, $N_2$ and $N_3$, but $\E[N_3/N^3\mid N>0]$ is in general not $\E N_3/\bar N^3$, and relating the two requires further information about the joint law of $N$ and $N_3$, such as concentration of $N$ around $\bar N$.
With thinning, the expected numbers of kept words, kept close pairs and kept far bad triples are $p\,\E N$, at most $p^2\E N_2$ and at most $p^3\E N_3$.
Only these first moments enter the proof, and no information on the fluctuations of $N$ is needed.

For Theorem~\ref{lower:main} we take $n=4m$ and $M=g^n$.
Then
\[
	M-M^2b_2-M^3b_3=M(1-Mb_2-M^2b_3),
\]
so Lemma~\ref{lower:two-stage} gives a trifferent code of blocklength $n$ with $(1-o(1))g^n$ words as soon as
\[
 \bar N\ge g^n,\qquad b_2=o(g^{-n}),\qquad b_3=o(g^{-2n})
\]
as $m\to\infty$.
Sections~\ref{lower:views}--\ref{lower:distance} construct a random set $S$ with these three properties, as the set of words accepted by a random system of overlapping constraints.

\subsection{An arbitrary system of overlapping views}
\label{lower:views}

Lemma~\ref{prelim:alteration} turns a law on short blocks into codes by repeating it independently on disjoint blocks.
A natural choice of block law is the uniform law on the blocks over $\F_9$ in which the labels occur with prescribed frequencies.
Whether such a law is useful depends not only on the number of allowed blocks but also on how many pairs of them are close and how many triples of them are bad.
We do not use disjoint blocks.
Instead, several constraints read each variable, so that every variable is subject to more than one restriction.

We build the code from words over $\F_9$ that satisfy many local constraints at once.
Recall from Section~\ref{lower:sampling} that the elements of $\F_9$ are called labels.
Each constraint reads $\ell$ of the $m$ positions of a word, which we call \emph{variables}, and tests the labels it reads against a fixed set $Q\subseteq\F_9^\ell$.
Different constraints may read the same variable, and this overlap is what makes the coordinates of an accepted word dependent.
Each reading also passes through its own permutation of $\F_9$.
The reason is explained after the definition.

\begin{definition}\label{lower:system}
Let $m$, $r$ and $\ell$ be positive integers.
\begin{enumerate}
\item An \emph{incidence structure} with $m$ variables and $r$ constraints of arity $\ell$ is a map
\[
 \Gamma\colon[r]\times[\ell]\to[m],\qquad (a,t)\mapsto j_{a,t}.
\]
The $r\ell$ pairs $(a,t)$ are its \emph{ports}.
The \emph{degree} of a variable $j\in[m]$ is the number of ports $(a,t)$ with $j_{a,t}=j$.
\item A \emph{local rule} is a set $Q\subseteq\F_9^\ell$.
\item A family of \emph{port permutations} consists of an invertible affine map $\pi_{a,t}(z)=\alpha z+\beta$ of $\F_9$, with $\alpha\ne\word{0}$, for every port $(a,t)$.
\end{enumerate}
For a word $x\in\F_9^m$, the \emph{view} of constraint $a$ is $(x_{j_{a,1}},\ldots,x_{j_{a,\ell}})\in\F_9^\ell$.
The word $x$ \emph{satisfies} constraint $a$ if
\[
 \bigl(\pi_{a,1}(x_{j_{a,1}}),\ldots,\pi_{a,\ell}(x_{j_{a,\ell}})\bigr)\in Q,
\]
and $x$ is \emph{accepted} if it satisfies every constraint $a\in[r]$.
\end{definition}

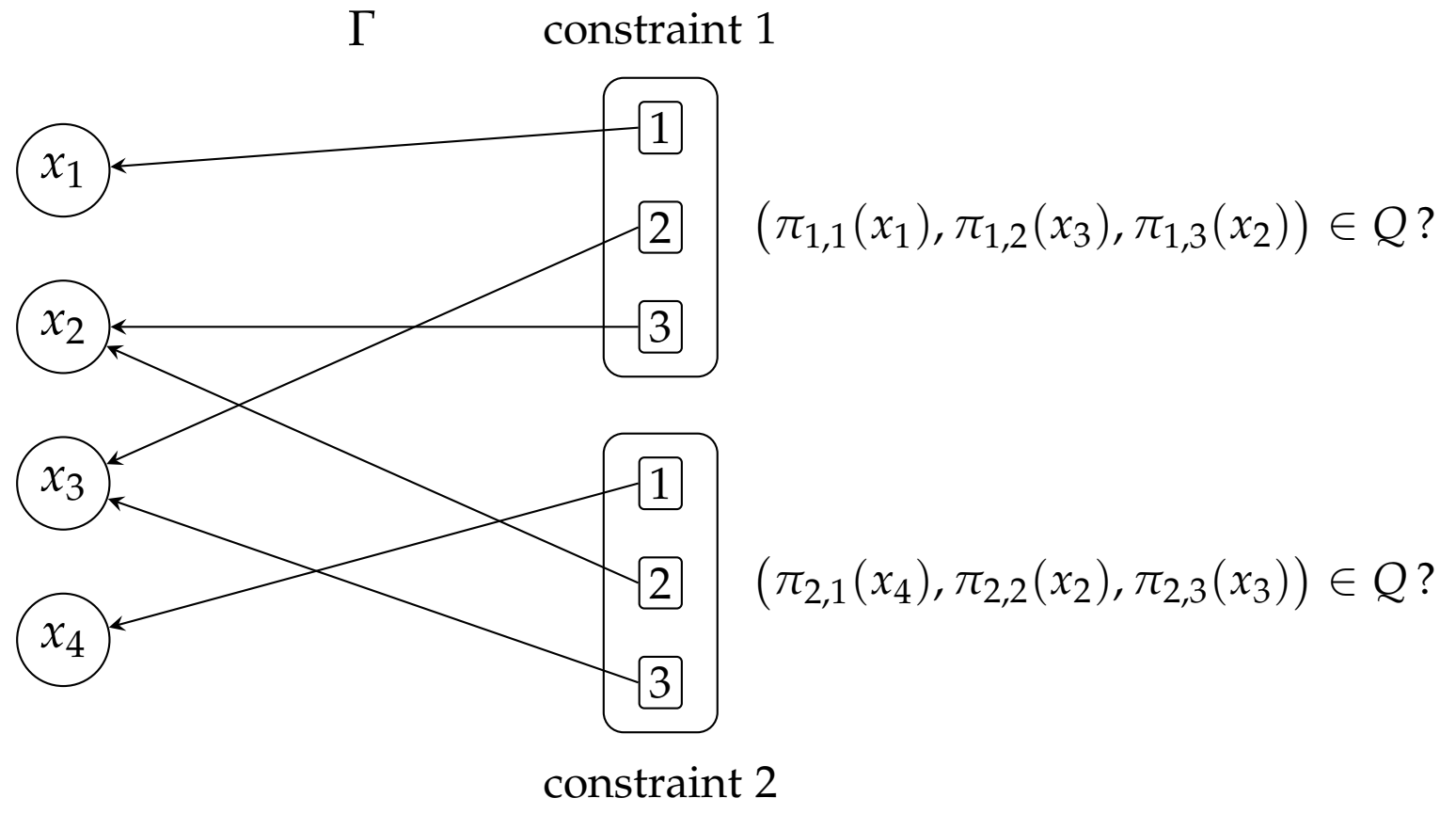
\begin{figure}[ht]
\centering
\begin{tikzpicture}[x=1cm,y=1cm,every node/.style={font=\small},
 var/.style={circle,draw,minimum size=6.5mm,inner sep=0pt},
 port/.style={draw,rounded corners=1pt,inner sep=2pt,font=\footnotesize},
 wire/.style={->,>=stealth,thin}]
\node[var] (x1) at (0,0) {$x_1$};
\node[var] (x2) at (0,-1.1) {$x_2$};
\node[var] (x3) at (0,-2.2) {$x_3$};
\node[var] (x4) at (0,-3.3) {$x_4$};
\node[port] (q11) at (4.2,0.3) {$1$};
\node[port] (q12) at (4.2,-0.4) {$2$};
\node[port] (q13) at (4.2,-1.1) {$3$};
\node[port] (q21) at (4.2,-2.2) {$1$};
\node[port] (q22) at (4.2,-2.9) {$2$};
\node[port] (q23) at (4.2,-3.6) {$3$};
\draw[rounded corners] (3.8,0.65) rectangle (4.6,-1.45);
\draw[rounded corners] (3.8,-1.85) rectangle (4.6,-3.95);
\node[anchor=west,font=\footnotesize] at (4.7,-0.4) {$\bigl(\pi_{1,1}(x_1),\pi_{1,2}(x_3),\pi_{1,3}(x_2)\bigr)\in Q\,?$};
\node[anchor=west,font=\footnotesize] at (4.7,-2.9) {$\bigl(\pi_{2,1}(x_4),\pi_{2,2}(x_2),\pi_{2,3}(x_3)\bigr)\in Q\,?$};
\node[font=\footnotesize] at (4.2,1.0) {constraint $1$};
\node[font=\footnotesize] at (4.2,-4.3) {constraint $2$};
\foreach \x/\q in {x1/q11, x3/q12, x2/q13, x4/q21, x2/q22, x3/q23}
 \draw[wire] (\q.west) -- (\x);
\node at (2.1,1.0) {$\Gamma$};
\end{tikzpicture}
\caption{An incidence structure with $m=4$ variables and $r=2$ constraints of arity $\ell=3$.
The arrows form the map $\Gamma$: the arrow from port $(a,t)$ points to the variable $j_{a,t}$ that the port reads.
The permutation $\pi_{a,t}$ belongs to the port.
Constraint $a$ applies it to the label of variable $j_{a,t}$ and tests the transformed view against $Q$, as shown on the right.
The variables $x_2$ and $x_3$ are read by both constraints, through different permutations.
The certified local rule has $\ell=14$.}
\label{lower:system-figure}
\end{figure}

A word $x\in\F_9^m$ gives the label $x_j$ to variable $j$.
The port $(a,t)$ is the $t$-th place at which constraint $a$ reads, and it reads the variable $j_{a,t}$.
So the view of constraint $a$ lists the labels of the variables it reads, in the order of its ports.
The constraint transforms the label at port $(a,t)$ by $\pi_{a,t}$ and tests the transformed view against $Q$.
In terms of maps, regard a word $x\in\F_9^m$ as the map $[m]\to\F_9$, $j\mapsto x_j$, and an element of $\F_9^\ell$ likewise as a map $[\ell]\to\F_9$.
Write $\Gamma_a\colon[\ell]\to[m]$ for the map $t\mapsto j_{a,t}$.
Then the view of constraint $a$ is the composite $x\circ\Gamma_a$, as in the left diagram below.
Let $\mathrm{AGL}(1,9)$ be the group of the $72$ invertible affine maps of $\F_9$, and write $\pi_a\colon[\ell]\to\mathrm{AGL}(1,9)$ for the map $t\mapsto\pi_{a,t}$.
Let $\mathrm{pr}_1$ and $\mathrm{pr}_2$ be the projections of $\mathrm{AGL}(1,9)\times\F_9$, and let $\mathrm{ev}(\pi,z)=\pi(z)$ be the evaluation map.
The map $(\pi_a,x\circ\Gamma_a)$ is the unique map that makes the right diagram commute.
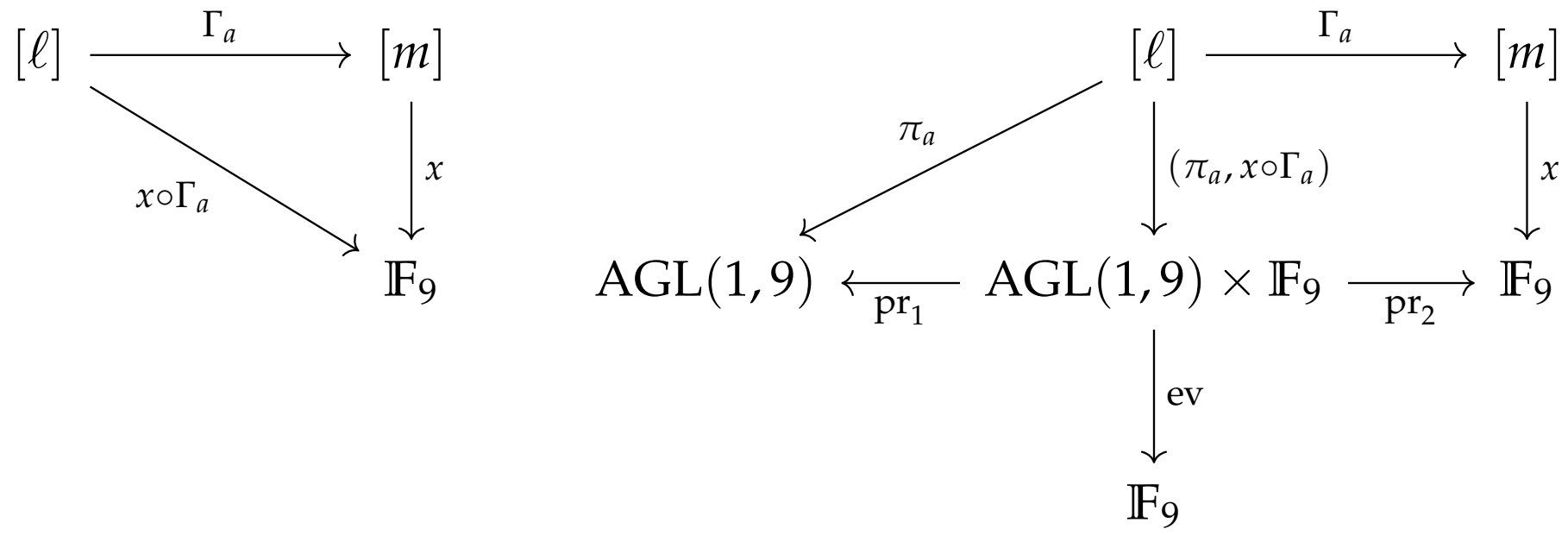
\[
\begin{tikzcd}[row sep=large,baseline=(\tikzcdmatrixname-1-1.base)]
 {[\ell]} \arrow[r,"\Gamma_a"] \arrow[dr,"x\circ\Gamma_a"'] & {[m]} \arrow[d,"x"] \\
 \phantom{\mathrm{AGL}(1,9)\times\F_9} & \F_9
\end{tikzcd}
\qquad
\begin{tikzcd}[row sep=large,baseline=(\tikzcdmatrixname-1-2.base)]
 & {[\ell]} \arrow[dl,"\pi_a"'] \arrow[d,"{(\pi_a,\,x\circ\Gamma_a)}"] \arrow[r,"\Gamma_a"] & {[m]} \arrow[d,"x"] \\
 \mathrm{AGL}(1,9) & \mathrm{AGL}(1,9)\times\F_9 \arrow[l,"\mathrm{pr}_1"] \arrow[r,"\mathrm{pr}_2"'] \arrow[d,"\mathrm{ev}"] & \F_9 \\
 & \F_9 &
\end{tikzcd}
\]
The transformed view is the composite $\mathrm{ev}\circ(\pi_a,x\circ\Gamma_a)$ down the middle of the right diagram, that is, $t\mapsto\pi_{a,t}\bigl(x(\Gamma_a(t))\bigr)$.
The map $\Gamma$ need not be injective.
A variable of degree $d$ is read at $d$ ports, by different constraints or several times by the same constraint, and its label is transformed separately at each of these ports.
Thus the permutations belong to the ports, not to the variables.
Figure~\ref{lower:system-figure} shows a small example.
Equivalently, $\Gamma$ is a hypergraph on the vertex set $[m]$ with $r$ labelled edges, in which each edge is an ordered list of $\ell$ vertices, possibly with repetitions.

The local rule $Q$ will be fixed.
We call a pair $\theta=(\Gamma,\pi)$, consisting of an incidence structure and a family of port permutations, a \emph{system}, and we write $\mathcal S_\theta\subseteq\F_9^m$ for its set of accepted words.
A \emph{random system} is a system chosen according to some probability distribution on systems.
Its set $\mathcal S_\theta$ of accepted words is then a random subset of $\F_9^m$.

We now explain the role of the port permutations.
Without them, a word is accepted if and only if every view lies in $Q$, and the number of accepted words depends on the detailed structure of $\Gamma$ and $Q$.
If $\Gamma$ is injective, the constraints read disjoint sets of variables, and exactly $|Q|^r9^{m-r\ell}$ words are accepted.
When a variable is read at several ports, however, one label enters several tests, and the tests are dependent.
Choosing the port permutations at random removes this dependence on average (Lemma~\ref{lower:ports-uniform}).

\begin{lemma}\label{lower:ports-uniform}
Fix an incidence structure $\Gamma$ and a local rule $Q$.
Choose the port permutations independently, each uniform among the $72$ invertible affine maps of $\F_9$.
Then, for every word $x\in\F_9^m$, the $r\ell$ transformed labels $\pi_{a,t}(x_{j_{a,t}})$ are independent and uniform on $\F_9$, and $x$ is accepted with probability $(|Q|/9^\ell)^r$.
Consequently, the number $N$ of accepted words satisfies
\begin{equation}\label{lower:mean-general}
 \bar N=\E N=9^m(|Q|/9^\ell)^r.
\end{equation}
\end{lemma}

\begin{proof}
Fix a word $x$.
First consider a single port $(a,t)$, and let $z=x_{j_{a,t}}$ be the label it reads.
For every $u\in\F_9$ and every $\alpha\ne\word{0}$, exactly one $\beta$ satisfies $\alpha z+\beta=u$.
So $\pi_{a,t}(z)=u$ for exactly $8$ of the $72$ affine maps, and $\pi_{a,t}(z)$ is uniform on $\F_9$.

Next consider all ports together.
The transformed label at port $(a,t)$ depends only on the permutation $\pi_{a,t}$, since $x$ is fixed.
The permutations at different ports are independent.
Hence the $r\ell$ transformed labels are independent.
This holds even when several ports read the same variable: they read the same label $z$, but transform it by independent permutations.

Finally, the transformed view of constraint $a$ consists of the $\ell$ transformed labels at its ports, so it is uniform on $\F_9^\ell$.
It lies in $Q$ with probability $|Q|/9^\ell$.
Different constraints have disjoint sets of ports, so their transformed views are independent, and $x$ is accepted with probability $(|Q|/9^\ell)^r$.
Summing over the $9^m$ words $x$ gives~\eqref{lower:mean-general}.
\end{proof}

In Lemma~\ref{lower:ports-uniform} the incidence structure is fixed, and the expectation is over the port permutations only.
So randomizing $\Gamma$ is not needed for $\bar N$.
For pairs and triples of words, the port permutations alone do not suffice to compute the expected counts, for the following reason.
For two words $x$ and $y$, the transformed labels at a port $(a,t)$ are $\pi_{a,t}(x_{j_{a,t}})$ and $\pi_{a,t}(y_{j_{a,t}})$.
They are equal if and only if $x_{j_{a,t}}=y_{j_{a,t}}$, because $\pi_{a,t}$ is a bijection.
So the port permutations do not make the transformed views of $x$ and $y$ independent.
The pattern of equal and unequal labels is carried from the variables to the ports by $\Gamma$, and Section~\ref{lower:counting} randomizes $\Gamma$ to control it.

From now on, $S=\mathcal S_\theta$ in Lemma~\ref{lower:two-stage}, for a random system $\theta=(\Gamma,\pi)$ whose port permutations are independent and uniform, as in Lemma~\ref{lower:ports-uniform}, and independent of $\Gamma$.
The incidence structure $\Gamma$ will also be random (Section~\ref{lower:certified}).
Then $\bar N$ is given by~\eqref{lower:mean-general}, and it remains to make $b_2$ and $b_3$ small.

\begin{remark}\label{lower:variable-permutations}
The permutations must be attached to the ports, not to the variables.
Suppose instead that there are permutations $\sigma_1,\ldots,\sigma_m$ of $\F_9$ with $\pi_{a,t}=\sigma_{j_{a,t}}$ for every port $(a,t)$.
Define $\phi\colon\F_9^m\to\F_9^m$ by $\phi(x)=(\sigma_1(x_1),\ldots,\sigma_m(x_m))$.
For every word $x$ and every $t\in[\ell]$,
\[
 \pi_{a,t}\bigl(x(\Gamma_a(t))\bigr)=\sigma_{\Gamma_a(t)}\bigl(x(\Gamma_a(t))\bigr)=\phi(x)\bigl(\Gamma_a(t)\bigr),
\]
so the transformed view of constraint $a$ in $x$ is $\phi(x)\circ\Gamma_a$, the view of constraint $a$ in $\phi(x)$.
Hence $x$ is accepted by the system $(\Gamma,\pi)$ if and only if $\phi(x)$ is accepted by the system in which every port permutation is the identity.
The map $\phi$ is a bijection.
It applies the same permutation to the labels of all words at a given position, so it preserves the label distance $d$ and maps bad triples to bad triples (Section~\ref{lower:sampling}).
Therefore the two systems have the same numbers of accepted words, of close pairs and of far bad triples.
Permutations attached to the variables thus have no effect on the construction.
Only the relation between the permutations at different ports of the same variable matters.
\end{remark}

\subsection{The statistics that the local rule must control}
\label{lower:counting}

At a position $i\in[\ell]$, two words $x,y\in\F_9^\ell$ are in state $\mathsf E$ if $x_i=y_i$ and in state $\mathsf N$ if $x_i\ne y_i$.
Write $\Sigma_2=\{\mathsf E,\mathsf N\}$.
For triples we record more than agreement.
Let $x,y,z\in\F_9^\ell$.
At a position $i\in[\ell]$ the labels $x_i,y_i,z_i$ are either all distinct, or they are in exactly one of the following four \emph{states}.
\begin{center}
\begin{tabular}{cl}\toprule
 State & Condition at position $i$\\\midrule
 $\mathsf A$ & $x_i=y_i=z_i$\\
 $\mathsf B_{12}$ & $x_i=y_i\ne z_i$\\
 $\mathsf B_{13}$ & $x_i=z_i\ne y_i$\\
 $\mathsf B_{23}$ & $y_i=z_i\ne x_i$\\\bottomrule
\end{tabular}
\end{center}
We call $\mathsf B_{12},\mathsf B_{13},\mathsf B_{23}$ the \emph{pair states}, and write $\Sigma_3=\{\mathsf A,\mathsf B_{12},\mathsf B_{13},\mathsf B_{23}\}$.
The triple $(x,y,z)$ is bad if no position has three distinct labels (Section~\ref{lower:sampling}), that is, if every position is in one of the four states.

The three pair states are kept separate because they give different pair distances.
For $s\in\Sigma_2$, let $c_s$ be the probability that two independent uniform labels are in state $s$, and for $s\in\Sigma_3$, let $c_s$ be the probability that three independent uniform labels are in state $s$.
So $c_{\mathsf E}=1/9$ and $c_{\mathsf N}=8/9$ are the probabilities that two independent uniform labels are equal or distinct.
Likewise, $c_{\mathsf A}=1/81$ is the probability that three independent uniform labels are all equal, and $c_{\mathsf B_{ik}}=8/81$ is the probability that the $i$-th and $k$-th are equal and the third is different.
We collect these numbers in the vectors $c_2=(c_{\mathsf E},c_{\mathsf N})=(1,8)/9$ and $c_3=(c_{\mathsf A},c_{\mathsf B_{12}},c_{\mathsf B_{13}},c_{\mathsf B_{23}})=(1,8,8,8)/81$.
The entries of $c_3$ sum to $25/81$.
The remaining probability $56/81$ is that of three distinct labels, which do not occur in a bad triple.

The \emph{pattern} of two words of $\F_9^\ell$ is the sequence $\sigma\in\Sigma_2^\ell$ of their states at the positions $1,\ldots,\ell$.
Likewise, the pattern of a bad triple of words of $\F_9^\ell$ is the sequence $\sigma\in\Sigma_3^\ell$ of its states.
Let $X,Y,Z$ be independent uniform elements of the local rule $Q\subseteq\F_9^\ell$.
For $\sigma\in\Sigma_2^\ell$, let $\nu_2(\sigma)$ be the probability that $(X,Y)$ has pattern $\sigma$.
For $\sigma\in\Sigma_3^\ell$, let $\nu_3(\sigma)$ be the probability that $(X,Y,Z)$ is bad and has pattern $\sigma$.

The construction never samples from the uniform law on $Q$.
The probabilities $\nu_2$ and $\nu_3$ enter through the expected counts of Section~\ref{lower:alteration}, as follows.
Fix two words of $\F_9^m$ and a constraint $a$, and let $\sigma\in\Sigma_2^\ell$ be the pattern of the views of the two words at $a$.
Choose the port permutations independently and uniformly, as in Lemma~\ref{lower:ports-uniform}.
Each $\pi_{a,t}$ is a bijection, so the transformed views of the two words also have pattern $\sigma$.
Position $t$ of the view of constraint $a$ is read through the port $(a,t)$, and we say that the port $(a,t)$ is in state $\sigma_t$.
Thus the port $(a,t)$ is in state $\mathsf E$ if the two words have the same label at the variable $j_{a,t}$, and in state $\mathsf N$ otherwise.
At a port in state $\mathsf E$, the two transformed labels are equal, and their common value is uniform on $\F_9$, as in the proof of Lemma~\ref{lower:ports-uniform}.
At a port in state $\mathsf N$, the pair of transformed labels is uniform among the $72$ ordered pairs of distinct labels, since for any two such pairs exactly one invertible affine map of $\F_9$ sends the first to the second.
The permutations at different ports are independent, so the pair of transformed views is uniform on the set of pairs in $(\F_9^\ell)^2$ with pattern $\sigma$.
This set has $\prod_{t=1}^\ell 81c_{\sigma_t}=81^\ell\prod_{t=1}^\ell c_{\sigma_t}$ elements, and $|Q|^2\nu_2(\sigma)$ of them lie in $Q^2$.
Hence both words satisfy constraint $a$ with probability
\begin{equation}
	\label{lower:local-acceptance}
	\Bigl(\frac{|Q|}{9^\ell}\Bigr)^{2}\frac{\nu_2(\sigma)}{\prod_{t=1}^\ell c_{\sigma_t}}.
\end{equation}
The same argument applies to three words whose views at $a$ form a bad triple with pattern $\sigma\in\Sigma_3^\ell$.
Each state $s\in\Sigma_3$ involves at most two distinct labels, and exactly $729c_s$ triples of labels are in state $s$.
So all three words satisfy constraint $a$ with probability
\begin{equation}\label{lower:local-acceptance-triple}
	\Bigl(\frac{|Q|}{9^\ell}\Bigr)^{3}\frac{\nu_3(\sigma)}{\prod_{t=1}^\ell c_{\sigma_t}}.
\end{equation}
This is why invertible affine maps suffice as port permutations.
Every state involves at most two distinct labels, and the affine maps of $\F_9$ act sharply transitively on ordered pairs of distinct labels.
The full symmetric group of $\F_9$ would give the same probabilities.
Apart from its size, the local rule thus enters the expected counts only through $\nu_2$ and $\nu_3$.

The patterns at different constraints are not independent.
For two words of $\F_9^m$, say that a variable $j\in[m]$ is in state $\mathsf E$ if the two words have the same label at $j$, and in state $\mathsf N$ otherwise.
For a bad triple of words of $\F_9^m$, the state of a variable in $\Sigma_3$ is defined in the same way.
Thus the state of a port $(a,t)$, as defined above, is the state of the variable $j_{a,t}$ that it reads, so all ports of a variable are in the same state.
For a fixed incidence structure, the expected counts therefore depend on how $\Gamma$ distributes the ports of each variable among the constraints, and computing them would require the detailed structure of $\Gamma$.
We avoid this by choosing $\Gamma$ at random as well (Section~\ref{lower:certified}).
The expected counts then depend on the local rule only through the sums of $\nu_2(\sigma)$ and $\nu_3(\sigma)$ over the patterns with prescribed numbers of positions in each state (Lemma~\ref{lower:coefficient-lemma}).
Section~\ref{lower:certified} collects these sums as the coefficients of two polynomials.
This random choice of $\Gamma$ is independent of the port permutations and of the thinning.

\subsection{A certified choice of local rule and wiring}
\label{lower:certified}

We use
\begin{equation}\label{lower:Q}
	Q=\{x\in\F_9^{14}:x\text{ has multiplicities }(4,3,1^7)\}.
\end{equation}
So the arity is $\ell=14$.
In each element of $Q$, one label occurs four times, a second label three times, and each of the remaining seven labels once.
Which two labels occur four and three times is not prescribed, so $|Q|=9\cdot8\cdot14!/(4!\,3!)$ (Appendix~\ref{lower:local-counts}).
The construction uses $Q$ only through its size and the probabilities $\nu_2$ and $\nu_3$ of Section~\ref{lower:counting}, and we chose $Q$ for these values.
Used on disjoint blocks, $Q$ is weaker.
Three independent uniform elements of $Q$ form a bad triple with probability approximately $5.99\cdot10^{-8}$.
Lemma~\ref{prelim:alteration}, applied to the image under $\Psi$ of the uniform law on $Q$, therefore gives only the base $1.16007\ldots$, below $g$.
Neither value is used in the proof.
The gain over this base comes from constraints that share variables.

Take $r=6k$ constraints, with $D=r\ell=84k$ ports in total, and $m=31k$ variables.
Put $m_2=9k$ and $m_3=22k$, so that $m=m_2+m_3$ and $D=2m_2+3m_3$.
The incidence structure $\Gamma$ is uniform among the maps $[r]\times[\ell]\to[m]$ under which the variables $1,\ldots,m_2$ have degree two and the variables $m_2+1,\ldots,m$ have degree three.
Equivalently, the ports are partitioned uniformly at random into labelled groups, one for each variable, of size two for the variables $1,\ldots,m_2$ and of size three for the others, and the ports in each group are the ports of its variable.
The method allows other degrees, and Remark~\ref{lower:uniform-wiring} explains why we prescribe degrees two and three rather than leave the degrees to a uniformly random map.
The encoded blocklength will be $n=4m=124k$.
The random system $\theta$ consists of this random incidence structure and of port permutations chosen independently and uniformly, as in Lemma~\ref{lower:ports-uniform}, independently of the incidence structure.
As before, $\mathcal S_\theta$ is its set of accepted words.
By~\eqref{lower:mean-general}, the expected number of accepted words is
\begin{equation}\label{lower:mean}
	\bar N=9^m\bigl(|Q|/9^{14}\bigr)^{6k}.
\end{equation}
Recall from Section~\ref{lower:alteration} that $b_2=\E N_2/\bar N^2$ and $b_3=\E N_3/\bar N^3$, where $N_2$ is the number of ordered pairs in $S^2$ at distance less than $\delta m$, and $N_3$ is the number of far bad ordered triples in $S^3$.
By the end of that section, Theorem~\ref{lower:main} follows once $\bar N\ge g^n$, $b_2=o(g^{-n})$ and $b_3=o(g^{-2n})$ as $m\to\infty$.
The next proposition gives these bounds for $S=\mathcal S_\theta$ and $\delta=3/8$, since the factors $(D+1)^2\rho_2^k$ and $(D+1)^6\rho_3^k$ tend to zero as $k\to\infty$.
\begin{proposition}[Certified counting estimates]\label{lower:budgets-prop}
Let $g=1.162579$ be the constant of Theorem~\ref{lower:main}, let $\delta=3/8$, and let $b_2$ and $b_3$ be the ratios~\eqref{lower:annealed-counts} for $S=\mathcal S_\theta$.
There are constants $\rho_2,\rho_3<1$, independent of $k$, such that for every $k\ge1$
\begin{equation}\label{lower:budgets}
	\bar N\ge g^n,\qquad
	b_2\le(D+1)^2\rho_2^k g^{-n},\qquad
	b_3\le(D+1)^6\rho_3^k g^{-2n}.
\end{equation}
\end{proposition}

The exponentially small factors $\rho_2^k$ and $\rho_3^k$ absorb the polynomial factors $(D+1)^2$ and $(D+1)^6$, so these are the three estimates required in Section~\ref{lower:alteration}.
The proof of Proposition~\ref{lower:budgets-prop} occupies the rest of this section and Section~\ref{lower:distance}.

We now derive the formula for the expected counts to which the certificates apply.
This is where the overlap of the constraints enters.
Let $\hat N_2=|\mathcal S_\theta|^2$ be the number of ordered pairs in $\mathcal S_\theta^2$, and let $\hat N_3$ be the number of ordered bad triples in $\mathcal S_\theta^3$.
So $\hat N_2=N^2$, where $N=|\mathcal S_\theta|$ as in Section~\ref{lower:alteration}.
The counts $N_2$ and $N_3$ of Section~\ref{lower:alteration} are the same counts with the distance conditions imposed.
Of the pairs counted by $\hat N_2$, the count $N_2$ keeps only the close ones, and of the triples counted by $\hat N_3$, the count $N_3$ keeps only the far ones.
Hence $N_2\le\hat N_2$ and $N_3\le\hat N_3$.
We first compute the expectations of $\hat N_2$ and $\hat N_3$ (Lemmas~\ref{lower:per-assignment-lemma} and~\ref{lower:coefficient-lemma}), and then obtain $N_2$ and $N_3$ by restricting the resulting sum.

Fix a pair ($j=2$) or a triple ($j=3$) of words, and record at each port the state of the labels it reads.
This port pattern can be described in two ways.
The constraints see it locally.
A constraint $a$ has the $\ell$ ports $(a,1),\ldots,(a,\ell)$, and their states form a pattern $\sigma\in\Sigma_j^\ell$.
By~\eqref{lower:local-acceptance} and~\eqref{lower:local-acceptance-triple}, the $j$ views at the constraint $a$ are all accepted with probability
\[
	\Bigl(\frac{|Q|}{9^\ell}\Bigr)^{j}\frac{\nu_j(\sigma)}{\prod_{t=1}^\ell c_{\sigma_t}}.
\]
The factor $\nu_j(\sigma)$ is where the local rule enters.
The denominator is a product over the $\ell$ ports of the constraint $a$, with the factor $c_{\sigma_t}$ for the port $(a,t)$.
Multiplying over the $r$ constraints therefore gives one factor $1/c_s$ for each of the $D=r\ell$ ports, where $s$ is the state of the port.

The variables see the pattern globally.
All ports of a variable read the labels at that variable, so they are all in the state of that variable.
Consider a variable of degree $d$ in state $s$.
Its $j$ labels can be chosen in $9^jc_s$ ways, so the variable contributes the factor $c_s$ to the number of $j$-tuples of words, apart from the factor $9^j$ common to all variables.
Its state is seen at its $d$ ports, so it contributes the factor $c_s^{-d}$ to the product of the acceptance probabilities over the constraints.
Its net weight is therefore $c_s^{1-d}$.

The wiring $\Gamma$ links the two descriptions.
The states of the variables determine how many ports are in each state, since a variable of degree $d$ in state $s$ puts $d$ ports in state $s$.
Which ports these are depends on $\Gamma$.
Since $\Gamma$ is uniform, all arrangements of the states over the $D$ ports with the given numbers $L_s$ of ports in each state $s$ are equally likely.
The expected count is therefore a sum over the vectors $L=(L_s)_{s\in\Sigma_j}$ of the product of three factors.
The first is the sum, over the arrangements with counts $L$, of the product of the factors $\nu_j(\sigma)$ over the constraints.
The second is the sum, over the assignments of states to the variables that produce the counts $L$, of the product of the net weights of the variables.
The third is one over the number of arrangements with counts $L$.
The two sums are coefficients of two polynomials, which we now define.

For the constraints, we use the \emph{local polynomials} of $Q$.
We introduce six indeterminates $w_s$, one for each state $s\in\Sigma_2\cup\Sigma_3$, and define
\begin{equation}\label{lower:local-polynomials}
\begin{aligned}
 H_Q(w_{\mathsf E},w_{\mathsf N})&=\sum_{\sigma\in\Sigma_2^\ell}\nu_2(\sigma)\prod_{t=1}^\ell w_{\sigma_t},\\
 F_Q(w_{\mathsf A},w_{\mathsf B_{12}},w_{\mathsf B_{13}},w_{\mathsf B_{23}})&=\sum_{\sigma\in\Sigma_3^\ell}\nu_3(\sigma)\prod_{t=1}^\ell w_{\sigma_t}.
\end{aligned}
\end{equation}
For $x,y,z\in\F_9^\ell$ and $s\in\Sigma_3$, let $N_s(x,y,z)$ be the number of positions in state $s$.
Grouping the outcomes of $(X,Y)$ and $(X,Y,Z)$ by their patterns shows that, equivalently, with $X,Y,Z$ independent and uniform on $Q$ as in Section~\ref{lower:counting},
\begin{align*}
 H_Q(w_{\mathsf E},w_{\mathsf N})&=\E\bigl[w_{\mathsf E}^{\#\{i:X_i=Y_i\}}w_{\mathsf N}^{\#\{i:X_i\ne Y_i\}}\bigr],\\
 F_Q(w_{\mathsf A},w_{\mathsf B_{12}},w_{\mathsf B_{13}},w_{\mathsf B_{23}})&=\E\Bigl[\mathbf1\{(X,Y,Z)\text{ is bad}\}\prod_{s\in\Sigma_3}w_s^{N_s(X,Y,Z)}\Bigr].
\end{align*}
Both are homogeneous polynomials of degree $\ell$ with nonnegative coefficients.
The coefficient of $\prod_sw_s^{\eta_s}$ in $H_Q$ or $F_Q$ is the sum of $\nu_2(\sigma)$ or $\nu_3(\sigma)$ over the patterns $\sigma$ with exactly $\eta_s$ positions in state $s$, for each $s$.
For example, the coefficient of $w_{\mathsf E}^iw_{\mathsf N}^{\ell-i}$ in $H_Q$ is the probability that $X$ and $Y$ agree at exactly $i$ positions, and $F_Q(1,1,1,1)$ is the probability that $(X,Y,Z)$ is bad.

Let $j=2$ for pairs and $j=3$ for bad triples, let $A=H_Q$ or $A=F_Q$ accordingly, write $w=(w_s)_{s\in\Sigma_j}$ for the indeterminates of $A$, and let $c=c_j$.
So
\[
 A(w)=\begin{cases}
  H_Q(w_{\mathsf E},w_{\mathsf N}) & \text{if } j=2,\\
  F_Q(w_{\mathsf A},w_{\mathsf B_{12}},w_{\mathsf B_{13}},w_{\mathsf B_{23}}) & \text{if } j=3.
 \end{cases}
\]
Let $v=(v_s)_{s\in\Sigma_j}$ be a second set of indeterminates.
For a vector $L=(L_s)_{s\in\Sigma_j}$ of nonnegative integers, write
\[
 |L|=\sum_{s\in\Sigma_j}L_s,\qquad
 w^L=\prod_{s\in\Sigma_j}w_s^{L_s},\qquad
 v^L=\prod_{s\in\Sigma_j}v_s^{L_s}.
\]
For a polynomial $P$, we write $[w^L]P(w)$ for the coefficient of the monomial $w^L$ in $P(w)$, and likewise $[v^L]P(v)$.
For an incidence structure with $r$ constraints in which $m_d$ variables have degree $d$, for each $d\ge1$, define
\begin{equation}\label{lower:VG}
	V(w)=A(w)^{r},\qquad
	R_d(v)=\sum_{s\in\Sigma_j} c_s^{1-d}v_s^d,\qquad
	G(v)=\prod_{d\ge1}R_d(v)^{m_d}.
\end{equation}
For the certified parameters, $V(w)=A(w)^{6k}$ and $G(v)=R_2(v)^{9k}R_3(v)^{22k}$.

\leavevmode\phantomsection\label{lower:degree-setting}%
We now fix the random system.
Let $Q\subseteq\F_9^\ell$ be a nonempty local rule, let $r\ge1$, put $D=r\ell$, and let $d_1,\ldots,d_m\ge1$ be integers with $d_1+\cdots+d_m=D$.
Let $m_d$ be the number of $i\in[m]$ with $d_i=d$, and let $V$ and $G$ be as in~\eqref{lower:VG}.
Let $\Gamma$ be uniform among the maps $[r]\times[\ell]\to[m]$ under which variable $i$ has degree $d_i$ for every $i\in[m]$.
The incidence structure of Section~\ref{lower:certified} is the case $d_i=2$ for $i\le m_2$ and $d_i=3$ for $i>m_2$.
Choose the port permutations independently and uniformly, as in Lemma~\ref{lower:ports-uniform}, and independently of $\Gamma$.
Let $\hat N_2$ and $\hat N_3$ be as above, for the accepted words of this random system, and recall that $\bar N=9^m(|Q|/9^\ell)^r$ by~\eqref{lower:mean-general}.
For $j\in\{2,3\}$ and an assignment $\omega\colon[m]\to\Sigma_j$ of states to the variables, let $\hat N_j(\omega)$ be the number of ordered $j$-tuples of accepted words whose state at every variable $i$ is $\omega(i)$.
Every ordered pair of words has exactly one assignment $\omega\colon[m]\to\Sigma_2$.
Every bad triple of words has exactly one assignment $\omega\colon[m]\to\Sigma_3$, and conversely a triple whose state at every variable lies in $\Sigma_3$ is bad, since two of its three labels are equal at every variable.
Hence
\[
 \hat N_j=\sum_{\omega\colon[m]\to\Sigma_j}\hat N_j(\omega)\qquad\text{for }j=2,3.
\]
Let $L(\omega)=(L_s(\omega))_{s\in\Sigma_j}$, where
\[
 L_s(\omega)=\sum_{i\in[m]:\,\omega(i)=s}d_i
\]
is the sum of the prescribed degrees of the variables in state $s$.
Every incidence structure in the support of $\Gamma$ gives variable $i$ the degree $d_i$, so $L(\omega)$ depends only on $\omega$ and the prescribed degrees, not on $\Gamma$.
A variable of degree $d_i$ is read at $d_i$ ports, all in its state, so $L_s(\omega)$ is the number of ports in state $s$, and $|L(\omega)|=D$.
We call $L(\omega)$ the \emph{port-count vector} of $\omega$.

\begin{lemma}[Per-assignment identity]\label{lower:per-assignment-lemma}
In this setting, let $j\in\{2,3\}$.
For every assignment $\omega\colon[m]\to\Sigma_j$,
\begin{equation}\label{lower:per-assignment}
 \frac{\E\hat N_j(\omega)}{\bar N^j}=[w^{L(\omega)}]V(w)\,\frac{\prod_{s\in\Sigma_j}L_s(\omega)!}{D!}\,\prod_{i=1}^mc_{\omega(i)}^{1-d_i}.
\end{equation}
\end{lemma}

\begin{proof}
First fix $\Gamma$.
The expectation is then over the port permutations only.
Fix an assignment $\omega\colon[m]\to\Sigma_j$.
Since $c_s$ is the probability that $j$ independent uniform labels are in state $s$, exactly $9^jc_s$ ordered $j$-tuples of labels are in state $s$.
A $j$-tuple of words of $\F_9^m$ amounts to a choice of one $j$-tuple of labels at each variable, and it has assignment $\omega$ if and only if the $j$-tuple of labels at each variable $i$ is in state $\omega(i)$.
So the number of $j$-tuples of words of $\F_9^m$ with assignment $\omega$ is
\[
	\prod_{i=1}^m9^jc_{\omega(i)}=9^{jm}\prod_{i=1}^mc_{\omega(i)}.
\]
Each of them is counted by $\hat N_j(\omega)$ if and only if its $j$ words are accepted.

Fix one of these $j$-tuples of words.
At constraint $a$, the port $(a,t)$ reads the variable $\Gamma_a(t)$, which is in state $\omega(\Gamma_a(t))$.
So the views of the $j$ words at $a$ have pattern
\[
	\sigma_a=\omega\circ\Gamma_a=(\omega(\Gamma_a(1)),\ldots,\omega(\Gamma_a(\ell))).
\]
By~\eqref{lower:local-acceptance} if $j=2$ and~\eqref{lower:local-acceptance-triple} if $j=3$, the $j$ words satisfy constraint $a$ with probability
\[
	\Bigl(\frac{|Q|}{9^\ell}\Bigr)^{j}\frac{\nu_j(\sigma_a)}{\prod_{t=1}^\ell c_{\omega(\Gamma_a(t))}}.
\]
The event that the $j$ words satisfy constraint $a$ depends only on the port permutations at the ports of $a$.
Different constraints have disjoint sets of ports, and the port permutations are independent, so these events are independent for different constraints.
Multiplying over the $r$ constraints, we find that the $j$ words are accepted with probability
\[
	\prod_{a=1}^r\Bigl(\frac{|Q|}{9^\ell}\Bigr)^{j}\frac{\nu_j(\sigma_a)}{\prod_{t=1}^\ell c_{\omega(\Gamma_a(t))}}
	=\Bigl(\frac{|Q|}{9^\ell}\Bigr)^{jr}\frac{\prod_{a=1}^r\nu_j(\sigma_a)}{\prod_{a=1}^r\prod_{t=1}^\ell c_{\omega(\Gamma_a(t))}}.
\]
The double product in the denominator runs over all ports, and each port contributes the factor of the state of the variable it reads.
Since variable $i$ is read by exactly $d_i$ ports, grouping the factors by variable gives
\[
 \prod_{a=1}^r\prod_{t=1}^\ell c_{\omega(\Gamma_a(t))}=\prod_{i=1}^m c_{\omega(i)}^{d_i}.
\]
So the $j$ words are accepted with probability
\[
	\Bigl(\frac{|Q|}{9^\ell}\Bigr)^{jr}\frac{\prod_{a=1}^r\nu_j(\sigma_a)}{\prod_{i=1}^m c_{\omega(i)}^{d_i}}.
\]
This probability depends on the $j$-tuple only through $\omega$.
By linearity of expectation, $\E[\hat N_j(\omega)\mid\Gamma]$ is therefore the number of $j$-tuples of words with assignment $\omega$ times this probability, that is,
\[
 \E[\hat N_j(\omega)\mid\Gamma]
 =9^{jm}\Bigl(\prod_{i=1}^m c_{\omega(i)}\Bigr)\Bigl(\frac{|Q|}{9^\ell}\Bigr)^{jr}\frac{\prod_{a=1}^r\nu_j(\sigma_a)}{\prod_{i=1}^m c_{\omega(i)}^{d_i}}.
\]
Dividing by $\bar N^j=9^{jm}(|Q|/9^\ell)^{jr}$, which follows from $\bar N=9^m(|Q|/9^\ell)^r$, we obtain
\begin{equation}\label{lower:fixed-wiring}
 \frac{\E[\hat N_j(\omega)\mid\Gamma]}{\bar N^j}
 =\Bigl(\prod_{i=1}^m c_{\omega(i)}^{1-d_i}\Bigr)\Bigl(\prod_{a=1}^r\nu_j(\omega\circ\Gamma_a)\Bigr).
\end{equation}
So a variable of degree $d$ in state $s$ carries the net weight $c_s^{1-d}$ described before the lemma.

Now let $\Gamma$ be random.
In~\eqref{lower:fixed-wiring}, the weight $\prod_{i=1}^m c_{\omega(i)}^{1-d_i}$ depends only on $\omega$ and the degrees, so only the factor $\prod_{a=1}^r\nu_j(\omega\circ\Gamma_a)$ is random.
Consider the tuple $(\omega\circ\Gamma_1,\ldots,\omega\circ\Gamma_r)$ of patterns.
Each pattern $\omega\circ\Gamma_a\in\Sigma_j^\ell$ has the $\ell$ entries $(\omega\circ\Gamma_a)(1),\ldots,(\omega\circ\Gamma_a)(\ell)$, so the tuple has $r\ell=D$ entries $(\omega\circ\Gamma_a)(t)$, one for each port $(a,t)$.
The entry $(\omega\circ\Gamma_a)(t)=\omega(\Gamma_a(t))$ is the state of the variable that the port $(a,t)$ reads, that is, the state of the port.
Whatever $\Gamma$ is, exactly $L_s(\omega)$ of these $D$ entries equal $s$, since $L_s(\omega)$ is the sum of the degrees of the variables in state $s$.
Composing $\Gamma$ with a permutation of the ports preserves the degrees, so it does not change the law of $\Gamma$.
Its effect on the tuple is to permute the $D$ entries, and these permutations act transitively on the tuples $(\sigma_1,\ldots,\sigma_r)\in(\Sigma_j^\ell)^r$ with exactly $L_s(\omega)$ entries $\sigma_a(t)$ equal to $s$ for each $s\in\Sigma_j$.
Hence $(\omega\circ\Gamma_1,\ldots,\omega\circ\Gamma_r)$ is uniform among these $D!/\prod_{s\in\Sigma_j}L_s(\omega)!$ tuples, and
\[
 \E\Bigl[\prod_{a=1}^r\nu_j(\omega\circ\Gamma_a)\Bigr]=\frac{\prod_{s\in\Sigma_j}L_s(\omega)!}{D!}\Biggl(\sum_{\substack{(\sigma_1,\ldots,\sigma_r)\in(\Sigma_j^\ell)^r\\\#\{(a,t):\,\sigma_a(t)=s\}=L_s(\omega)\text{ for all }s\in\Sigma_j}}\ \prod_{a=1}^r\nu_j(\sigma_a)\Biggr).
\]
On the other hand, by~\eqref{lower:local-polynomials},
\[
 V(w)=A(w)^r=\prod_{a=1}^r\Biggl(\sum_{\sigma_a\in\Sigma_j^\ell}\nu_j(\sigma_a)\prod_{t=1}^\ell w_{\sigma_a(t)}\Biggr).
\]
Expanding the product gives one term for each tuple $(\sigma_1,\ldots,\sigma_r)\in(\Sigma_j^\ell)^r$, namely $\prod_{a=1}^r\nu_j(\sigma_a)$ times the monomial $\prod_{s\in\Sigma_j}w_s^{\#\{(a,t):\,\sigma_a(t)=s\}}$.
So the coefficient $[w^{L(\omega)}]V(w)$ is the sum of $\prod_{a=1}^r\nu_j(\sigma_a)$ over the tuples with $\#\{(a,t):\sigma_a(t)=s\}=L_s(\omega)$ for all $s\in\Sigma_j$, which is the sum above, and taking the expectation of~\eqref{lower:fixed-wiring} over $\Gamma$ gives~\eqref{lower:per-assignment}.
\end{proof}

\begin{lemma}[Coefficient identity]\label{lower:coefficient-lemma}
In the setting of Lemma~\ref{lower:per-assignment-lemma}, let $j\in\{2,3\}$.
Then
\begin{equation}\label{lower:coefficient-sum}
 \frac{\E\hat N_j}{\bar N^j}=
 \sum_{\substack{L\colon\Sigma_j\to\mathbb Z_{\ge0}\\|L|=D}}[w^L]V(w)\,[v^L]G(v)\,
       \frac{\prod_{s\in\Sigma_j}L_s!}{D!}.
\end{equation}
\end{lemma}

\begin{proof}
Recall that $\hat N_j$ is the sum of $\hat N_j(\omega)$ over all assignments $\omega\colon[m]\to\Sigma_j$.
Multiplying out $G(v)=\prod_{i=1}^mR_{d_i}(v)$ amounts to choosing a state $\omega(i)$ for each variable $i$, with weight $c_{\omega(i)}^{1-d_i}$ and monomial $v_{\omega(i)}^{d_i}$.
The choice $\omega$ therefore contributes $\prod_{i=1}^mc_{\omega(i)}^{1-d_i}$ times the monomial $\prod_{i=1}^mv_{\omega(i)}^{d_i}=v^{L(\omega)}$.
Hence, for every $L\colon\Sigma_j\to\mathbb Z_{\ge0}$ with $|L|=D$,
\[
 [v^L]G(v)=\sum_{\substack{\omega\colon[m]\to\Sigma_j\\L(\omega)=L}}\ \prod_{i=1}^mc_{\omega(i)}^{1-d_i}.
\]
Summing~\eqref{lower:per-assignment} over all $\omega$ and grouping the assignments by $L(\omega)$, we obtain
\[
 \frac{\E\hat N_j}{\bar N^j}
 =\sum_{\substack{L\colon\Sigma_j\to\mathbb Z_{\ge0}\\|L|=D}}[w^L]V(w)\,\frac{\prod_{s\in\Sigma_j}L_s!}{D!}\Biggl(\sum_{\substack{\omega\colon[m]\to\Sigma_j\\L(\omega)=L}}\ \prod_{i=1}^mc_{\omega(i)}^{1-d_i}\Biggr),
\]
and the inner sum is $[v^L]G(v)$, which gives~\eqref{lower:coefficient-sum}.
\end{proof}

The factorial ratio in~\eqref{lower:coefficient-sum} is the probability that the uniformly random port pattern with counts $L$ is any one prescribed pattern.
For a general local rule $Q$, it is the uniform choice of $\Gamma$ that reduces the sum to total state counts.
No symmetry of $Q$ under permutations of positions is needed.

\subsection{Keeping distance in the count}
\label{lower:distance}

Recall from Section~\ref{lower:certified} that $\hat N_j$ is the sum of $\hat N_j(\omega)$ over all assignments $\omega\colon[m]\to\Sigma_j$.
The counts $N_2$ and $N_3$ of Section~\ref{lower:alteration} are the same sums, restricted to the assignments of close pairs and of far bad triples, which we now describe.
For a pair $(x,y)$ of words of $\F_9^m$ with assignment $\omega\colon[m]\to\Sigma_2$, the words differ exactly at the variables in state $\mathsf N$, so $d(x,y)=|\omega^{-1}(\mathsf N)|$.
Let $\Omega_2$ be the set of assignments $\omega\colon[m]\to\Sigma_2$ with $|\omega^{-1}(\mathsf N)|<\delta m$.
Thus a pair is close if and only if its assignment lies in $\Omega_2$.

For a bad triple $(x,y,z)$ of words of $\F_9^m$ with assignment $\omega\colon[m]\to\Sigma_3$, the words $x$ and $y$ differ exactly at the variables in state $\mathsf B_{13}$ or $\mathsf B_{23}$, and similarly for the other two pairs.
So
\[
 d(x,y)=|\omega^{-1}\{\mathsf B_{13},\mathsf B_{23}\}|,\qquad
 d(x,z)=|\omega^{-1}\{\mathsf B_{12},\mathsf B_{23}\}|,\qquad
 d(y,z)=|\omega^{-1}\{\mathsf B_{12},\mathsf B_{13}\}|.
\]
Let $\Omega_3$ be the set of assignments $\omega\colon[m]\to\Sigma_3$ for which each of these three numbers is at least $\delta m$.
Thus a bad triple is far if and only if its assignment lies in $\Omega_3$.
Consequently $N_j=\sum_{\omega\in\Omega_j}\hat N_j(\omega)$ for $j=2,3$.

A distance between words counts variables, whereas the port-count vector $L(\omega)$ counts ports, and a variable of degree $d$ accounts for $d$ ports.
So whether $\omega$ lies in $\Omega_j$ is not determined by $L(\omega)$, and the restriction to $\Omega_j$ cannot be imposed on the coefficients of~\eqref{lower:coefficient-sum} directly.
We therefore bound the indicator of $\Omega_j$ by a product of weights, one for each variable.
Since $G$ is a product of one factor for each variable, these weights can be built into its factors, and this gives the polynomial $G_t$ below.

The weights depend on a parameter $t$, which we call the \emph{tilt}.
For $j=2$, the tilt is a number $t$ with $0<t\le1$.
For $j=3$, the tilt is a triple $t=(t_{12},t_{13},t_{23})$ of numbers with $t_{12},t_{13},t_{23}\ge1$, one for each pair of positions in an ordered triple of words.
For a triple $(x,y,z)$, the number $t_{12}$ goes with the pair $(x,y)$, $t_{13}$ with $(x,z)$, and $t_{23}$ with $(y,z)$.
The tilt is a parameter, not a function of the words.
The same tilt is used for all pairs, or all triples, of words, and in Lemma~\ref{lower:weighted-bound} it is chosen separately for each port-count vector.
Put
\[
 \Lambda(t)=\log t\quad\text{if }j=2,\qquad
 \Lambda(t)=\log t_{12}+\log t_{13}+\log t_{23}\quad\text{if }j=3.
\]
To each state $s\in\Sigma_j$ we assign the \emph{weight} $\tau_s$ in the following tables.
\begin{center}
\begin{tabular}{ccc}\toprule
 $s$ & $\mathsf E$ & $\mathsf N$\\\midrule
 $\tau_s$ & $1$ & $t$\\\bottomrule
\end{tabular}
\qquad\qquad
\begin{tabular}{ccccc}\toprule
 $s$ & $\mathsf A$ & $\mathsf B_{12}$ & $\mathsf B_{13}$ & $\mathsf B_{23}$\\\midrule
 $\tau_s$ & $1$ & $t_{13}t_{23}$ & $t_{12}t_{23}$ & $t_{12}t_{13}$\\\bottomrule
\end{tabular}
\end{center}
In words, $\tau_s$ is the product of the tilts of the pairs of words that differ in state $s$.
For example, in state $\mathsf B_{12}$ the first two words agree and the third differs from both, so exactly the pairs $(x,z)$ and $(y,z)$ differ, and $\tau_{\mathsf B_{12}}=t_{13}t_{23}$.
For $j=2$ there is a single pair, which differs in state $\mathsf N$ and agrees in state $\mathsf E$.
The weights enter the variable polynomials as follows:
\begin{equation}\label{lower:weighted-G}
 R_{d,t}(v)=\sum_{s\in\Sigma_j}c_s^{1-d}\tau_sv_s^d,\qquad
 G_t(v)=\prod_{d\ge1}R_{d,t}(v)^{m_d}.
\end{equation}
Each term of $R_{d,t}$ contains the weight $\tau_s$ once, not once for each of the $d$ ports, because a variable contributes at most one to each distance.
If all tilts equal $1$, then all weights equal $1$, and $R_{d,t}$ and $G_t$ reduce to $R_d$ and $G$.

Now let $(x,y)$ be a pair, or $(x,y,z)$ a bad triple, of words of $\F_9^m$.
Its assignment is the map $\omega\colon[m]\to\Sigma_j$ that sends each variable $i$ to the state of the labels $x_i,y_i$, or $x_i,y_i,z_i$, at $i$.
Give each variable $i$ the weight $\tau_{\omega(i)}$ of its state.
The tilt of a pair of words is a factor of $\tau_{\omega(i)}$ exactly when that pair differs at the variable $i$.
In the product of the weights over all variables, it therefore occurs once for each variable at which the pair differs, that is,
\[
 \prod_{i=1}^m\tau_{\omega(i)}=t^{d(x,y)}\quad\text{if }j=2,\qquad
 \prod_{i=1}^m\tau_{\omega(i)}=t_{12}^{d(x,y)}\,t_{13}^{d(x,z)}\,t_{23}^{d(y,z)}\quad\text{if }j=3.
\]
If $\omega\in\Omega_2$, then $d(x,y)<\delta m$, and since $t\le1$,
\[
 \exp(-\delta m\Lambda(t))\prod_{i=1}^m\tau_{\omega(i)}=t^{d(x,y)-\delta m}\ge1.
\]
If $\omega\in\Omega_3$, then each of the three distances is at least $\delta m$, and since $t_{12},t_{13},t_{23}\ge1$,
\[
 \exp(-\delta m\Lambda(t))\prod_{i=1}^m\tau_{\omega(i)}=t_{12}^{d(x,y)-\delta m}\,t_{13}^{d(x,z)-\delta m}\,t_{23}^{d(y,z)-\delta m}\ge1.
\]
In both cases,
\begin{equation}\label{lower:weight-domination}
 \exp(-\delta m\Lambda(t))\prod_{i=1}^m\tau_{\omega(i)}\ge1\qquad\text{for every }\omega\in\Omega_j.
\end{equation}
Let $\mathcal L_j=\{L(\omega):\omega\in\Omega_j\}$ be the set of port-count vectors of the assignments in $\Omega_j$, that is, of the assignments of close pairs if $j=2$ and of far bad triples if $j=3$.

\begin{lemma}[Weighted coefficient bound]\label{lower:weighted-bound}
In the setting of Lemma~\ref{lower:per-assignment-lemma}, let $\delta>0$ and $j\in\{2,3\}$, and let $b_j$ be the ratio~\eqref{lower:annealed-counts} for the set of accepted words.
For each $L\in\mathcal L_j$, choose a tilt $t_L$ as above.
Then
\begin{equation}\label{lower:tilted-sum}
 b_j\le\sum_{L\in\mathcal L_j}\exp(-\delta m\Lambda(t_L))\,[w^L]V(w)\,[v^L]G_{t_L}(v)\,\frac{\prod_{s\in\Sigma_j}L_s!}{D!}.
\end{equation}
\end{lemma}

\begin{proof}
Since $N_j=\sum_{\omega\in\Omega_j}\hat N_j(\omega)$, identity~\eqref{lower:per-assignment} gives
\[
 b_j=\frac{\E N_j}{\bar N^j}=\sum_{L\in\mathcal L_j}[w^L]V(w)\,\frac{\prod_{s\in\Sigma_j}L_s!}{D!}
 \sum_{\substack{\omega\in\Omega_j\\L(\omega)=L}}\ \prod_{i=1}^mc_{\omega(i)}^{1-d_i}.
\]
Fix $L\in\mathcal L_j$ and write $t=t_L$.
By~\eqref{lower:weight-domination}, each term of the inner sum satisfies
\[
 \prod_{i=1}^mc_{\omega(i)}^{1-d_i}\le\exp(-\delta m\Lambda(t))\prod_{i=1}^mc_{\omega(i)}^{1-d_i}\tau_{\omega(i)}.
\]
Extending the inner sum to all assignments $\omega$ with $L(\omega)=L$ only adds nonnegative terms.
As in the proof of Lemma~\ref{lower:coefficient-lemma}, multiplying out $G_t(v)=\prod_{i=1}^mR_{d_i,t}(v)$ shows that the sum of $\prod_{i=1}^mc_{\omega(i)}^{1-d_i}\tau_{\omega(i)}$ over these assignments is $[v^L]G_t(v)$.
Hence the inner sum is at most $\exp(-\delta m\Lambda(t))[v^L]G_t(v)$, which proves~\eqref{lower:tilted-sum}.
\end{proof}

The tilt may depend on $L$ because~\eqref{lower:weight-domination} holds for every tilt, and in the proof it is applied separately to the inner sum of each $L$, before the sum over $L$ is taken.

It remains to bound the sum~\eqref{lower:tilted-sum} for large $D$.
We bound each coefficient by evaluating its polynomial at a positive point, and the factorial ratio by a multinomial estimate.
The resulting bound on the term of $L$ is exponential in $D$, with a rate $B_j$ that depends on the \emph{port profile} $h=L/D$, that is, on the vector of fractions $h_s=L_s/D$ of ports in each state $s$, on the evaluation points and on the tilt.
The task is then the following.
For every port profile $h$ of a close pair or a far bad triple, find positive evaluation points $w,v$ and a tilt $t$, a \emph{witness}, for which the rate is less than a threshold $\kappa$ fixed by the target base $g$.
The witness may vary with $h$.
Since the rate is convex in $h$, one witness serves for a whole polytope of port profiles if the rate is less than $\kappa$ at the vertices of the polytope, so finitely many polytopes and witnesses suffice.
Let $J=|\Sigma_j|$, so that $J=2$ if $j=2$ and $J=4$ if $j=3$, and put
\[
 \alpha_d=dm_d/D,\qquad e=m/D=\sum_{d\ge1}\alpha_d/d.
\]
So $\alpha_d$ is the fraction of ports that belong to variables of degree $d$, and $e$ is the number of variables per port.
For a port profile $h=(h_s)_{s\in\Sigma_j}$, by which we mean any vector with $h_s\ge0$ and $\sum_{s\in\Sigma_j}h_s=1$, for positive $w$ and $v$, and for a tilt $t$ as above, define
\begin{equation}\label{lower:cell-bound}
 B_j(h;w,v,t)=\frac1\ell\log A(w)
   +\sum_{d\ge1}\frac{\alpha_d}{d}\log R_{d,t}(v)
   -\delta e\Lambda(t)
   +\sum_{s\in\Sigma_j}h_s\log\frac{h_s}{w_sv_s},
\end{equation}
with $0\log0=0$.
The four terms of $B_j$ come from the four factors of the term of $L$ in~\eqref{lower:tilted-sum}.
The first comes from the constraint coefficient $[w^L]V(w)$, and the second from the variable coefficient $[v^L]G_t(v)$.
The third comes from the distance factor $\exp(-\delta m\Lambda(t))$.
The last combines the estimate of the factorial ratio with the factors $w^{-L}$ and $v^{-L}$ that the evaluation at $w$ and $v$ leaves over.
For a fixed witness $(w,v,t)$, the rate is continuous and convex in $h$ on the closed simplex, since $h_s\log h_s$ is convex and the remaining terms are affine in $h$.

\begin{lemma}[Finite cell bound]\label{lower:cell-lemma}
Fix $\ell$, a nonempty local rule $Q\subseteq\F_9^\ell$, $\delta>0$ and $j\in\{2,3\}$.
Fix also nonnegative numbers $\alpha_d$ for $d\ge1$, finitely many of them nonzero, with $\sum_{d\ge1}\alpha_d=1$.
Call $D$ \emph{admissible} if $D/\ell$ and all $\alpha_dD/d$ are integers.
For admissible $D$, consider the random system of Lemma~\ref{lower:per-assignment-lemma} with $r=D/\ell$ constraints and $m_d=\alpha_dD/d$ variables of degree $d$, and let $b_j$ be its ratio~\eqref{lower:annealed-counts} and $\mathcal L_j$ its set of port-count vectors.
Let $\kappa$ be a real number.
Suppose that we are given finitely many closed polytopes in the simplex of port profiles, the \emph{cells}, and for each cell a \emph{witness} $(w,v,t)$, consisting of $w,v\in(0,\infty)^{\Sigma_j}$ and a tilt $t$, such that:
\begin{enumerate}
\item for every admissible $D$ and every $L\in\mathcal L_j$, the port profile $L/D$ lies in some cell;
\item $B_j(h;w,v,t)<\kappa$ at every vertex $h$ of every cell, where $(w,v,t)$ is the witness of that cell.
\end{enumerate}
Then there is $\zeta>0$ such that, for every admissible $D$,
\begin{equation}\label{lower:finite-norm}
 b_j\le (D+1)^{2(J-1)}
           \exp\bigl(D(\kappa-\zeta)\bigr).
\end{equation}
\end{lemma}

\begin{proof}
Let $L\in\mathcal L_j$, put $h=L/D$, and let $(w,v,t)$ be the witness of a cell that contains $h$.
Use $t_L=t$ in~\eqref{lower:tilted-sum}.
Since $V$ and $G_t$ have nonnegative coefficients,
\[
 [w^L]V(w)\le A(w)^rw^{-L},\qquad
 [v^L]G_t(v)\le G_t(v)v^{-L}.
\]
We also use the multinomial estimate
\begin{equation}\label{lower:multinomial}
 \frac{\prod_{s\in\Sigma_j}L_s!}{D!}\le(D+1)^{J-1}
       \exp\Bigl(D\sum_{s\in\Sigma_j}h_s\log h_s\Bigr),
\end{equation}
which is proved below.
Since $r=D/\ell$, $m_d=\alpha_dD/d$, $m=eD$ and $L=Dh$,
\begin{align*}
 A(w)^r&=\exp\Bigl(\frac D\ell\log A(w)\Bigr),&
 G_t(v)&=\exp\Bigl(D\sum_{d\ge1}\frac{\alpha_d}{d}\log R_{d,t}(v)\Bigr),\\
 w^{-L}v^{-L}&=\exp\Bigl(-D\sum_{s\in\Sigma_j}h_s\log(w_sv_s)\Bigr),&
 \exp(-\delta m\Lambda(t))&=\exp(-D\delta e\Lambda(t)).
\end{align*}
Multiplying these bounds and comparing with~\eqref{lower:cell-bound}, we find that the term of $L$ in~\eqref{lower:tilted-sum} satisfies
\[
 \exp(-\delta m\Lambda(t))\,[w^L]V(w)\,[v^L]G_t(v)\,\frac{\prod_{s\in\Sigma_j}L_s!}{D!}
 \le(D+1)^{J-1}\exp\bigl(DB_j(h;w,v,t)\bigr).
\]

Every point $h$ of a cell is a convex combination $\sum_i\lambda_ih^{(i)}$ of the vertices $h^{(i)}$ of the cell, so convexity gives
\[
 B_j(h;w,v,t)\le\sum_i\lambda_iB_j(h^{(i)};w,v,t)\le\max_iB_j(h^{(i)};w,v,t)<\kappa,
\]
with the witness $(w,v,t)$ of that cell.
The vertex values do not depend on $D$, and there are finitely many of them, all less than $\kappa$.
So there is $\zeta>0$ such that $B_j\le\kappa-\zeta$ on every cell, with the witness of that cell.
Finally, there are at most $(D+1)^{J-1}$ vectors $L$ with $|L|=D$, so
\[
 b_j\le\sum_{\substack{L\colon\Sigma_j\to\mathbb Z_{\ge0}\\|L|=D}}(D+1)^{J-1}\exp\bigl(D(\kappa-\zeta)\bigr)
 \le(D+1)^{2(J-1)}\exp\bigl(D(\kappa-\zeta)\bigr),
\]
which is~\eqref{lower:finite-norm}.

It remains to prove~\eqref{lower:multinomial}.
Consider $D$ independent trials, each with outcome $s$ with probability $h_s$, and let $P$ be the probability that each outcome $s$ occurs exactly $L_s$ times.
Since $L=Dh$,
\[
 P=\frac{D!}{\prod_{s\in\Sigma_j}L_s!}\prod_{s\in\Sigma_j}h_s^{L_s}
  =\frac{D!}{\prod_{s\in\Sigma_j}L_s!}\exp\Bigl(D\sum_{s\in\Sigma_j}h_s\log h_s\Bigr),
\]
so~\eqref{lower:multinomial} says that $P\ge(D+1)^{-(J-1)}$.
Order $\Sigma_j$ as $s_1,\ldots,s_J$, and for $i=1,\ldots,J$ put
\[
 D_i=D-L_{s_1}-\cdots-L_{s_{i-1}}=D(h_{s_i}+\cdots+h_{s_J}),\qquad
 p_i=\frac{h_{s_i}}{h_{s_i}+\cdots+h_{s_J}}\quad\text{if }D_i>0.
\]
Revealing the counts of $s_1,s_2,\ldots$ one at a time gives
\[
 P=\prod_{\substack{1\le i\le J-1\\D_i>0}}\binom{D_i}{L_{s_i}}p_i^{L_{s_i}}(1-p_i)^{D_i-L_{s_i}}.
\]
The factor of $i$ is the probability that a binomial random variable with $D_i$ trials and success probability $p_i$ equals its mean $D_ip_i=L_{s_i}$, which is an integer.
A binomial random variable $X$ with $D'$ trials and integer mean $D'p$ takes the value $D'p$ with probability at least $(D'+1)^{-1}$.
Indeed, if $p\in\{0,1\}$ then $X=D'p$ surely.
Otherwise
\[
 \frac{\Pr[X=k+1]}{\Pr[X=k]}=\frac{(D'-k)p}{(k+1)(1-p)}
\]
is at least $1$ for $k<D'p$ and less than $1$ for $k\ge D'p$, so $D'p$ is a most likely value among the $D'+1$ possible values.
As $D_i\le D$, each factor is at least $(D+1)^{-1}$, and $P\ge(D+1)^{-(J-1)}$.
\end{proof}

For the certified parameters, $\ell=14$, $D=84k$, and
\[
 \alpha_2=3/14,\qquad\alpha_3=11/14,\qquad
 e=31/84,\qquad \delta=3/8.
\]
Since $n=4m=4eD$, the bounds on $b_2$ and $b_3$ in~\eqref{lower:budgets} have the form~\eqref{lower:finite-norm}, with $\rho_j=\exp(-84\zeta_j)<1$, where $\zeta_j$ is the constant of Lemma~\ref{lower:cell-lemma} for $j$, and
\begin{equation}\label{lower:certified-inequalities}
 \kappa_2=-4e\log g,\qquad \kappa_3=-8e\log g,
 \qquad g=1162579/1000000.
\end{equation}

\begin{proposition}[Certified vertex inequalities]\label{cert:lower-norms}
For the certified parameters and $j\in\{2,3\}$ the following hold.
Let the \emph{pair domain} be the closed segment of port profiles $(h_{\mathsf E},h_{\mathsf N})$, and let the \emph{triple domain} be the closed simplex of port profiles with $h_{\mathsf B_{12}}\ge h_{\mathsf B_{13}}\ge h_{\mathsf B_{23}}$.
Each of these domains is covered by finitely many rational simplices, the \emph{cells}.
Each cell is either excluded or retained.
A cell is excluded only if there is one condition in~\eqref{lower:pair-domain}, for $j=2$, or in~\eqref{lower:domain}, for $j=3$, that every vertex of the cell violates strictly.
Each retained cell carries a witness $(w,v,t)$ with rational entries, where $w,v\in(0,\infty)^{\Sigma_j}$ and $t$ is a tilt, that is, $0<t\le1$ if $j=2$ and $t_{12},t_{13},t_{23}\ge1$ if $j=3$.
With these witnesses, $B_j(h;w,v,t)<\kappa_j$ at every vertex $h$ of every retained cell.
\end{proposition}

Proposition~\ref{cert:lower-norms} is a finite computation in exact rational arithmetic with rigorous logarithm bounds.
Numerical optimization was used only to find the witnesses.
Appendix~\ref{cert:appendix} describes the data and how they are verified.

\begin{proof}[Proof of Proposition~\ref{lower:budgets-prop}]
We apply Lemma~\ref{lower:cell-lemma} with $\kappa=\kappa_j$, the admissible $D$ being the multiples of $84$.

For $j=2$, take as cells the retained cells of Proposition~\ref{cert:lower-norms}, with their witnesses.
Let $L\in\mathcal L_2$.
By Appendix~\ref{lower:domain-details}, the port profile $L/D$ satisfies~\eqref{lower:pair-domain}, so it lies in some cell of the cover.
That cell is retained, because in an excluded cell every vertex violates the same condition of~\eqref{lower:pair-domain} strictly, and hence by convexity so does every point.
This is hypothesis~1 of Lemma~\ref{lower:cell-lemma}, and hypothesis~2 is the vertex inequality of Proposition~\ref{cert:lower-norms}.

For $j=3$, the triple domain is only the part of the simplex where the pair-state coordinates are in decreasing order, and we use the symmetry under permutations of the three words.
Let $\pi$ be a permutation of the three positions of an ordered triple of words.
It fixes the state $\mathsf A$, and it permutes the three pairs of positions, and with them the pair states $\mathsf B_{12},\mathsf B_{13},\mathsf B_{23}$ and the tilts $t_{12},t_{13},t_{23}$, in the same way.
Write $\pi h$, $\pi w$, $\pi v$ and $\pi t$ for the vectors with permuted coordinates.
Then
\[
 B_3(\pi h;\pi w,\pi v,\pi t)=B_3(h;w,v,t).
\]
Indeed, $F_Q(\pi w)=F_Q(w)$, because permuting the three independent uniform elements $X,Y,Z$ of $Q$ does not change their joint law and permutes the pair states in the same way.
The constants $c_s$ agree on the three pair states.
Since $\tau_s$ is the product of the tilts of the pairs that differ in state $s$, and $\pi$ maps these pairs to the pairs that differ in state $\pi s$, the weight of $\pi s$ under $\pi t$ is the weight of $s$ under $t$, so $R_{d,\pi t}(\pi v)=R_{d,t}(v)$.
Moreover $\Lambda(\pi t)=\Lambda(t)$, and the last sum in~\eqref{lower:cell-bound} is only reordered.
Also, $\pi t$ is again a tilt.
Take as cells the images under all six permutations $\pi$ of the retained cells, with the correspondingly permuted witnesses.
By the displayed invariance, hypothesis~2 of Lemma~\ref{lower:cell-lemma} follows from Proposition~\ref{cert:lower-norms}.
For hypothesis~1, let $L\in\mathcal L_3$ and $h=L/D$.
Choose $\pi$ such that $\pi h$ lies in the triple domain.
The conditions~\eqref{lower:domain} are invariant under permuting the pair states, so $\pi h$ satisfies them by Appendix~\ref{lower:domain-details}.
As for $j=2$, $\pi h$ therefore lies in a retained cell $C$, and $h$ lies in the cell $\pi^{-1}C$.

Lemma~\ref{lower:cell-lemma} therefore gives the bounds on $b_2$ and $b_3$ in~\eqref{lower:budgets}.
For the first bound, \eqref{lower:mean} gives $\bar N=(|Q|^6/9^{53})^k$, and the exact integer comparison
\begin{equation}\label{lower:first-moment-check}
 |Q|^6\,1000000^{124}>9^{53}\,1162579^{124}
\end{equation}
gives $\bar N\ge g^{124k}=g^n$.
\end{proof}

\begin{proof}[Proof of Theorem~\ref{lower:main}]
By Proposition~\ref{lower:budgets-prop}, $\bar N\ge g^n$, so Lemma~\ref{lower:two-stage} applies with $\delta=3/8$ and $M=g^n$.
The other two bounds of the proposition give
\[
 M^2b_2+M^3b_3
 \le M\bigl((D+1)^2\rho_2^k+(D+1)^6\rho_3^k\bigr)=o(M).
\]
The lemma therefore gives a trifferent code of blocklength $n=4m$ with $(1-o(1))g^n$ words.
Its minimum distance is at least $3\delta m=9m/8=9n/32$, which is the by-product stated after the theorem.
For an arbitrary blocklength $n\ge124$, use the largest multiple of $124$ that is at most $n$, and pad with fixed coordinates.
This changes the size by only a constant factor.
Padding does not change distances, so the padded code has minimum distance at least $(9/32)(n-123)$, and the codes have asymptotic relative minimum distance at least $9/32$.
\end{proof}

The $o(1)$ term decays slowly.
If $\zeta_3$ is taken to be the certified margin $1.28\cdot10^{-8}$ of Appendix~\ref{cert:certificates}, the factor $(D+1)^6\rho_3^k$ falls below $1$ only for $k$ above about $1.3\cdot10^8$, that is, for $n$ above about $1.6\cdot10^{10}$.

\begin{remark}[Why the degrees are prescribed]\label{lower:uniform-wiring}
The most natural random incidence structure is a uniform map $\Gamma\colon[r]\times[\ell]\to[m]$, under which every port reads an independent uniform variable.
The incidence structure of Section~\ref{lower:certified} is this map conditioned on the degrees of the variables.
Without the conditioning the method gives much less.

Identity~\eqref{lower:fixed-wiring} holds for every fixed $\Gamma$, also when some variables have degree zero.
Fix $j\in\{2,3\}$ and an assignment $\omega\colon[m]\to\Sigma_j$, and let $\varphi_s$ be the fraction of the variables in state $s$.
Under the uniform map the states at the $D$ ports are independent, and each equals $s$ with probability $\varphi_s$.
The factor $\prod_{i=1}^mc_{\omega(i)}^{-d_i}$ in~\eqref{lower:fixed-wiring} is a product over the ports, so taking the expectation over $\Gamma$ gives
\[
 \frac{\E\hat N_j(\omega)}{\bar N^j}=\prod_{s\in\Sigma_j}c_s^{m\varphi_s}\,A\bigl((\varphi_s/c_s)_{s\in\Sigma_j}\bigr)^r,
\]
with $A=H_Q$ for $j=2$ and $A=F_Q$ for $j=3$.
Take $r/m=6/31$ and $\delta=3/8$, as in Section~\ref{lower:certified}, and let $m$ be divisible by $93$, so that $r=6m/31$ and $m/3$ are integers.
The $\binom{m}{m/3}$ assignments with $(\varphi_{\mathsf E},\varphi_{\mathsf N})=(2/3,1/3)$ are close, so
\[
 b_2\ge\binom{m}{m/3}\Bigl(\frac19\Bigr)^{2m/3}\Bigl(\frac89\Bigr)^{m/3}H_Q(6,3/8)^{6m/31}
 =\exp\bigl((-0.54986\ldots+o(1))m\bigr).
\]
Since $\log(9/5)=0.58778\ldots$, every $M\ge(9/5)^m$ satisfies $M^2b_2>M$ for large $m$, and the right-hand side of~\eqref{lower:joint-alteration} is negative.
As $n=4m$, Lemma~\ref{lower:two-stage} then does not reach even the base $(9/5)^{1/4}$ of K\"orner and Marton.
Numerical optimization gives the following approximate values for the largest $g$ with $\bar N\ge g^n$, $b_2=o(g^{-n})$ and $b_3=o(g^{-2n})$.
They are not certified.
For $r/m=6/31$ and the best $\delta$, it is about $1.1446$.
Over all $r/m$ and $\delta$, it is about $1.1609$, attained near $r/m=0.115$ and $\delta=0.385$.
This exceeds $(9/5)^{1/4}$ but falls short of $g=\gstar$.

The loss comes from the spread of the degrees.
Under the uniform map the degree of a variable is binomial with $D$ trials and success probability $1/m$, so the degrees are asymptotically Poisson with mean $D/m$.
At $D/m=84/31$, about $7\%$ of the variables are read by no constraint and about $29\%$ have degree at least four.
A variable of degree $d$ in state $s$ carries the net weight $c_s^{1-d}$, and
\[
 \Bigl(\frac{c_{\mathsf E}}{c_{\mathsf N}}\Bigr)^{1-d}=\Bigl(\frac{c_{\mathsf A}}{c_{\mathsf B_{12}}}\Bigr)^{1-d}=8^{d-1}.
\]
So the net weight of a variable of degree $d$ is $8^{d-1}$ times larger in the agreement state $\mathsf E$ or $\mathsf A$ than in the state $\mathsf N$ or in one of the states $\mathsf B_{12},\mathsf B_{13},\mathsf B_{23}$.
Two words that agree at such a variable agree at all $d$ of its ports, while their distance decreases by only one.
A variable of degree zero is read by no constraint, so its label is free.
Three independent uniform labels form a bad triple with probability $25/81$, so for a fixed incidence structure such a variable multiplies $\E\hat N_3/\bar N^3$, which counts all bad triples, by exactly $25/81$.
This is the factor per variable of the construction of K\"orner and Marton.
Since $n=4m$, the target $b_3=o(g^{-2n})$ asks for a factor below $g^{-8}\approx0.2997$ per variable on average, and $25/81\approx0.3086$ exceeds it.
So variables of degree zero must be compensated by the others.
This is a heuristic, not a bound, since $b_3$ counts only far bad triples.

Concentrated degrees are the best choice for the bound of Lemma~\ref{lower:cell-lemma}.
Lemmas~\ref{lower:per-assignment-lemma}--\ref{lower:cell-lemma} extend to degree zero with $R_{0,t}(v)=\sum_{s\in\Sigma_j}c_s\tau_s$.
The degree term of~\eqref{lower:cell-bound} is then $D^{-1}\sum_{i=1}^m\log R_{d_i,t}(v)$, where $d_1,\ldots,d_m$ are the degrees.
For fixed $v$ and $t$, the function
\[
 d\mapsto\log R_{d,t}(v)=\log\sum_{s\in\Sigma_j}c_s\tau_s\exp\bigl(d\log(v_s/c_s)\bigr)
\]
is convex on $[0,\infty)$, as the logarithm of a sum of exponentials of affine functions of $d$.
Moving one port from a variable of degree $d$ to a variable of degree at most $d-2$ therefore does not increase the degree term.
Among all degree sequences with $m$ variables and $D$ ports, the degree term is thus smallest when every degree is $\lfloor D/m\rfloor$ or $\lceil D/m\rceil$.
For $D/m=84/31$ these are the $9k$ variables of degree two and the $22k$ of degree three.
The other terms of~\eqref{lower:cell-bound}, and $\bar N$, depend on $m$ and $D$ but not on the degrees.
So, for every port profile $h$ and every witness, the degrees of Section~\ref{lower:certified} minimize $B_j(h;w,v,t)$ among all degree sequences with the same $m$ and $D$.
This compares the bounds of Lemma~\ref{lower:cell-lemma}, not the ratios $b_2$ and $b_3$ themselves, and the set of port profiles that the cells must cover also depends on the degrees.
The numerical values in this remark come from floating-point computations that the supplementary programs do not reproduce, and none of them is used in the proof of Theorem~\ref{lower:main}.
\end{remark}

\section{A half-power saving for arbitrary codes}
\label{sec:general}

A binary box contains at most two words of a trifferent code, and averaging this over all boxes gives the classical bound $T(n)\le2(3/2)^n$.
We also use the empty boxes.
Certain six-word configurations of codewords, the objects below, are each seen from exactly three empty boxes.
Two overlapping objects determine an empty box, and Lemma~\ref{gen:coherence} shows that all overlaps at a given object determine the same one.
Counting the empty boxes from which an object is seen without overlaps then bounds the number of boxes containing two codewords.

\begin{theorem}\label{thm:half-power}
Every ternary trifferent code of blocklength $n$ has size $O(n^{-1/2}(3/2)^n)$.
\end{theorem}

\subsection{Boxes, contacts, and a six-word configuration}

Throughout this section $n\ge1$, and $C\subseteq\F_3^n$ is a code of size $M$.
A \emph{mask} is a word $q\in\F_3^n$.
A word $x$ \emph{avoids} $q$ if $x_i\ne q_i$ for every $i\in[n]$, so the words avoiding $q$ form the binary box avoiding $q$.
Let $f(q)$ be the number of codewords that avoid $q$.
Three words \emph{triffer} if they have a rainbow coordinate.
Three words avoid a common mask if and only if they do not triffer, since a coordinate admits a symbol different from all three entries exactly when it is not rainbow.
So $C$ is trifferent if and only if $f(q)\le2$ for every mask $q$.
From now on $C$ is trifferent.
A mask $q$ is a \emph{hole} if $f(q)=0$, and a \emph{double mask} if $f(q)=2$.
The two codewords avoiding a double mask are its \emph{owners}.
Let $Z$ and $D$ be the numbers of holes and of double masks.

We first record anchored packing, the discrete specialization of Della Fiore's centered inequality \cite[Lemma~1]{DellaFiore2026}.
Fix an anchor $x\in C$.
For distinct $y,y'\in C\setminus\{x\}$, no mask is avoided by all of $x,y,y'$, so the sets of masks avoided by both $x$ and $y$, for $y\in C\setminus\{x\}$, are disjoint.
A mask avoided by $x$ and $y$ has two possible symbols at each coordinate where $x$ and $y$ agree and one at each coordinate where they differ, so there are $2^{n-d(x,y)}$ such masks.
All of them lie among the $2^n$ masks avoided by $x$.
Therefore
\begin{equation}\label{gen:packing}
 \sum_{y\in C\setminus\{x\}}2^{-d(x,y)}\le1.
\end{equation}
This holds for every trifferent code of every blocklength.
If equality holds for an anchor $x$, then no word $z\notin C$ can be added without destroying trifference.
Indeed, some mask is avoided by both $x$ and $z$, and by equality this mask is also avoided by a codeword $y\ne x$, so $x,y,z$ do not triffer.

Let $q$ be a mask and $i\in[n]$.
A codeword $x$ is a \emph{contact} of $q$ in direction $i$ if $x_i=q_i$ and $x_j\ne q_j$ for every $j\ne i$.
There are at most two contacts in each direction.
Indeed, three contacts in direction $i$ agree at $i$, and at every other coordinate $j$ they avoid $q_j$ and so take at most two values, so they do not triffer.
If $q$ has two contacts $u,v$ in direction $i$, then $q$ is a hole.
Indeed, a codeword $z$ avoiding $q$ would differ from $u$ and $v$ at $i$, and $u,v,z$ would not triffer.
At $i$ two of them equal $q_i$, and at every other coordinate all three avoid $q$.
At a hole $q$, let $J(q)$ be the set of directions in which $q$ has two contacts, and put $d_q=|J(q)|$.
The two contacts in a direction $i\in J(q)$ form the \emph{contact pair} in direction $i$.

A \emph{completed object}, or simply an \emph{object}, is a pair $(S,m)$, where $S\subseteq[n]$ has three elements and $m$ is a mask on the coordinates outside $S$, with the following property.
The \emph{restriction} of $(S,m)$ is the set of codewords $x$ with $x_i\ne m_i$ for every $i\notin S$.
We require that it consist of six codewords whose projections to $S$ form, up to a permutation of the coordinates of $S$ and independent permutations of the symbols in each coordinate, the set
\begin{equation}\label{gen:six}
 A_0=\word{000},\quad A_1=\word{111},\quad A_2=\word{222},
 \qquad B_0=\word{012},\quad B_1=\word{120},\quad B_2=\word{201}.
\end{equation}
We call $S$ the \emph{support} and $m$ the \emph{outside mask} of the object, and we label its six codewords as in~\eqref{gen:six}.
Let $H$ be the number of objects.

The restriction of any pair $(S,m)$ projects to a trifferent set.
Indeed, three of its words take at most two values at each coordinate outside $S$, so their rainbow coordinate lies in $S$.
When the restriction has at least three words, the projection is injective, since two words with the same projection and any third word would have no rainbow coordinate.
The five distances from each word of~\eqref{gen:six} to the others are $2,2,2,3,3$, so~\eqref{gen:packing} is an equality for this trifferent set of blocklength three.
So if six codewords avoid $m$ outside $S$ and project to a copy of~\eqref{gen:six}, the restriction contains no seventh codeword, whose projection would be a seventh word, and $(S,m)$ is an object.
The triples $A_0A_1A_2$ and $B_0B_1B_2$ are rainbow at all three coordinates.
They are the only such triples, since a triple with two words of one kind and one of the other is rainbow at exactly one coordinate.
We call them the two \emph{rainbow triples} of the object.

Each object $(S,m)$ has six \emph{centers}.
The center of one of its codewords $w$ is the mask that equals $w$ on $S$ and $m$ outside $S$.
The codewords avoiding a center lie in the restriction.
A direct check on~\eqref{gen:six} shows that exactly two of the six words avoid the center of $w$ on $S$, namely the other two words of the rainbow triple of $w$.
So the center is a double mask, and these two words are its owners.
The word $w$ itself equals the center on $S$, and we call it the \emph{candidate} of the center.
The remaining three words agree with the center on $S$ in exactly one coordinate each, one for each coordinate of $S$, so they are contacts of the center, one in each direction of $S$.

Let $\Omega$ be the graph on the set of objects in which two objects are adjacent when they have a common center and their supports meet in two coordinates.

\subsection{Three representations of each object}
\label{gen:reps}

Let $q$ be a hole.
For distinct $i,j\in J(q)$, the two words $u,v$ of the contact pair in direction $i$ satisfy $u_j\ne v_j$.
Indeed, let $w$ be a contact in direction $j$.
At $i$ both $u$ and $v$ equal $q_i$, and at every coordinate outside $\{i,j\}$ all three words avoid $q$.
So the rainbow coordinate of $u,v,w$ is $j$, where $w_j=q_j$, and therefore $u_j\ne v_j$.

At each coordinate $j\in J(q)$, encode the symbol $q_j$ by $0$ and the two other symbols by the signs $+1$ and $-1$, in some order.
Name the words of each contact pair $u_i^+$ and $u_i^-$, in some order.
By the previous paragraph there are signs $a_{ij}$, for distinct $i,j\in J(q)$, such that
\begin{equation}\label{gen:skew}
 u_i^\pm(i)=0,\qquad u_i^\pm(j)=\pm a_{ij},\qquad
 a_{ij}=-a_{ji}\in\{1,-1\},
\end{equation}
where the last condition holds after a suitable choice of the encodings and names, as we now show.
Fix a reference direction $i_0\in J(q)$.
For every $i\ne i_0$, name the pair in direction $i$ so that $a_{ii_0}=+1$, and choose the encoding at $i$ so that $a_{i_0i}=-1$.
Now let $i,j\ne i_0$ be distinct.
The words $u_i^+,u_j^+,u_{i_0}^+$ have the entries $+1,+1,0$ at $i_0$, are nonzero at the other coordinates of $J(q)\setminus\{i,j\}$, and avoid $q$ outside $J(q)$.
So their rainbow coordinate is $i$ or $j$.
At $i$ their entries are $0,a_{ji},-1$, and at $j$ they are $a_{ij},0,-1$.
Hence $a_{ji}=+1$ or $a_{ij}=+1$.
The same argument with $u_{i_0}^-$, whose entries at $i$ and $j$ are $+1$, gives $a_{ji}=-1$ or $a_{ij}=-1$.
Together, $a_{ij}=-a_{ji}$.

Any three directions $a,b,c\in J(q)$ give an object.
Put $S=\{a,b,c\}$ and let $m$ be the restriction of $q$ to the coordinates outside $S$.
The six contacts in the directions of $S$ avoid $m$ outside $S$, and by~\eqref{gen:skew} their projections to $S$ are determined by the three signs $a_{ab},a_{ac},a_{bc}$.
Exchanging the two nonzero symbols at a coordinate $i\in S$, together with the names $u_i^+$ and $u_i^-$, preserves~\eqref{gen:skew} and negates $a_{ij}$ and $a_{ji}$ for both $j\in S\setminus\{i\}$.
These operations change the number of minus signs among $a_{ab},a_{ac},a_{bc}$ by an even number and act transitively on the sign vectors of a given parity.
Exchanging two coordinates of $S$ changes the parity.
So all choices of the three signs give the same six-word set, up to permutations of coordinates and of symbols.
The six words of~\eqref{gen:six} arise in this way, as the contacts of the mask $\word{021}$ in the trifferent set~\eqref{gen:six} itself.
Hence the six contacts project to a copy of~\eqref{gen:six}, and $(S,m)$ is an object.

Conversely, let $(S,m)$ be an object.
A \emph{hole representation} of it is a hole $q$ that equals $m$ outside $S$ and satisfies $S\subseteq J(q)$.
For such $q$ the six contacts in the directions of $S$ avoid $m$ outside $S$, so they are the six words of the object, and each of them agrees with $q$ on $S$ in exactly one coordinate.
In the coordinates of~\eqref{gen:six}, the masks on $S$ with this property are exactly
\begin{equation}\label{gen:three-holes}
 \word{021},\qquad\word{102},\qquad\word{210}.
\end{equation}
Indeed, $A_0$, $A_1$ and $A_2$ each agree with such a mask exactly once, so it uses each symbol once, and the remaining permutation words $\word{012}$, $\word{120}$ and $\word{201}$ are excluded because they equal $B_0$, $B_1$ and $B_2$.
Conversely, each mask~\eqref{gen:three-holes}, extended by $m$, has the six words as contacts, two in each direction of $S$, so it is a hole $q$ with $S\subseteq J(q)$.
So every object has exactly three hole representations.
Every pair of a hole $q$ and a three-element subset of $J(q)$ gives an object of which $q$ is a hole representation, and therefore
\begin{equation}\label{gen:objects}
 \sum_{q\in\F_3^n:\,f(q)=0}\binom{d_q}{3}=3H.
\end{equation}

Fix a hole $q$.
Let $\Gamma_q$ be the graph whose vertices are the three-element subsets of $J(q)$, each standing for the object it gives, in which two vertices are adjacent when their objects are adjacent in $\Omega$.
Distinct vertices give distinct objects, since their supports differ.
Let $S=\{a,b,x\}$ and $T=\{a,b,y\}$ be vertices with $S\cap T=\{a,b\}$.
We claim that they are adjacent if and only if
\begin{equation}\label{gen:edge}
 a_{xa}a_{xb}=a_{ya}a_{yb}.
\end{equation}
Every center of the object of $S$ equals $q$ outside $S$, in particular at $y$, and every center of the object of $T$ equals $q$ at $x$.
So a common center equals $q$ outside $\{a,b\}$.
The centers of the object of $S$ that equal $q$ at $x$ are those of the two contacts in direction $x$, since the other four words differ from $q$ at $x$.
At $a$ and $b$ these two centers have the sign patterns $\pm(a_{xa},a_{xb})$, and elsewhere they equal $q$.
Likewise, the centers of the object of $T$ that equal $q$ at $y$ have the sign patterns $\pm(a_{ya},a_{yb})$ at $a$ and $b$.
Two pairs of opposite sign vectors $\pm(s_1,s_2)$ and $\pm(s_1',s_2')$ coincide exactly when $s_1s_2=s_1's_2'$, which proves~\eqref{gen:edge}.

\subsection{A conflict closes into a four-coordinate square}
\label{gen:square}

Take an edge of $\Omega$ between objects with supports $S$ and $T$, and let $r$ be a common center.
Put $U=S\cup T$, which has four coordinates.
The mask $r$ equals the outside mask of each object outside its support, so the two outside masks agree outside $U$.
Since $f(r)=2$, the owners of $r$ are the same for both objects.
They are two words of a rainbow triple of each object, so they differ at every coordinate of $S$ and of $T$, and both avoid $r$.
So at every coordinate of $U$ the mask $r$ and the two owners have three distinct symbols, and we may normalize $r$ to $\word{0000}$ and the owners to $\word{1111}$ and $\word{2222}$ on $U$.
Each object contributes the three contacts of $r$ in the directions of its support.
Since $r$ is not a hole, it has at most one contact in each direction, so the two objects have the same contacts in the two directions of $S\cap T$.
So $r$ has exactly one contact in each direction of $U$, giving four contacts in the directions of $U$.
Each object also contributes the candidate of $r$, which equals $\word{0}$ on its support and is nonzero at the coordinate of $U$ outside it.
These eight words are distinct, since the sets of coordinates of $U$ at which they equal $\word{0}$ are empty for the owners, singletons for the contacts, and the two distinct sets $S$ and $T$ for the candidates.
All eight avoid $r$ outside $U$, so they lie in the set of codewords that avoid $r$ at every coordinate outside $U$, which we call the restriction of the edge.
As before, its projection to $U$ is injective and trifferent.

Write $c_i$ for the contact in direction $i\in U$, and encode its entries $\word{1}$ and $\word{2}$ at the coordinates $j\ne i$ by $+1$ and $-1$.
For distinct $i,j\in U$, the words $c_i$, $c_j$ and $\word{1111}$ have no rainbow coordinate outside $\{i,j\}$, since outside $U$ all three avoid $r$ and on $U\setminus\{i,j\}$ all three entries are nonzero.
At $i$ their entries are $\word{0},c_j(i),\word{1}$, and at $j$ they are $c_i(j),\word{0},\word{1}$, so $c_j(i)=\word{2}$ or $c_i(j)=\word{2}$.
The owner $\word{2222}$ likewise gives $c_j(i)=\word{1}$ or $c_i(j)=\word{1}$.
Hence $c_i(j)=-c_j(i)$, and the contacts form a skew four-by-four sign table.
Orient the complete graph on $U$ by putting $i\to j$ when $c_i(j)=+1$.

Consider the candidate with support $S$, a contact $c_i$ with $i\in S$, and an owner.
At $i$ two of the three entries are $\word{0}$.
At the coordinate of $U$ outside $S$ all three entries are nonzero, and outside $U$ all three words avoid $r$.
So the rainbow coordinate lies in $S\setminus\{i\}$, where the entries are $\word{0}$, $c_i(j)$ and the owner's symbol.
With the owner $\word{1111}$ this gives $c_i(j)=\word{2}$ for some $j\in S\setminus\{i\}$, and with $\word{2222}$ it gives $c_i(j)=\word{1}$ for some $j\in S\setminus\{i\}$.
So every vertex of $S$ has exactly one out-neighbour in $S$, and $S$ is a cyclic triangle of the tournament.
The same holds for $T$.

Let $S\cap T=\{a,b\}$, and after exchanging the names of $a$ and $b$ assume that $a\to b$.
The cyclic triangles then force $b\to c\to a$ for the third vertex $c$ of $S$, and $b\to d\to a$ for the third vertex $d$ of $T$.
The edge between $c$ and $d$ has two orientations, which are exchanged by exchanging $c$ and $d$.
Up to these permutations of coordinates the tournament, and hence the contacts, are determined, and a representative is
\begin{equation}\label{gen:eight}
\begin{array}{llll}
 A=\word{0001},&B=\word{0020},&C=\word{0211},&D=\word{1022},\\
 E=\word{1111},&F=\word{2102},&G=\word{2110},&H=\word{2222}.
\end{array}
\end{equation}
Here $C,D,F,G$ are the contacts and $A,B$ the candidates.
A candidate equals $r$ on its support, so only its entry at the remaining coordinate of $U$ is free.
For $A=\word{000}v$, the triple $A,F,G$ can only be rainbow at the last coordinate, which forces $v=\word{1}$.
For $B=\word{00}w\word{0}$, the triple $B,F,G$ can only be rainbow at the third coordinate, which forces $w=\word{2}$.
Each word of~\eqref{gen:eight} has the distances $2,2,3,3,3,4,4$ to the other seven, so~\eqref{gen:packing} is an equality for this trifferent set of blocklength four.
As before, the restriction of the edge contains no ninth codeword.

Consider the objects whose supports lie in $U$ and whose outside masks equal $r$ outside $U$.
Such an object has support $U\setminus\{t\}$ for some $t\in U$, and its six codewords are the words of~\eqref{gen:eight} that avoid its outside mask at $t$.
Each column of~\eqref{gen:eight} contains one symbol twice and each of the other two symbols three times, so only the symbol used twice leaves six words.
This gives at most one object for each $t$, and a direct check shows that each of the four six-word sets projects to a copy of~\eqref{gen:six}.
The four objects are as follows, where the slash separates their two rainbow triples, and the fixed coordinate of an object is the coordinate $t$ of $U$ outside its support, with the value of its outside mask at $t$.
\begin{center}
\begin{tabular}{lll}\toprule
Object & Triples & Fixed coordinate and value\\\midrule
$P_0$ & $AGH/BCF$ & $0=\word{1}$\\
$P_1$ & $ADG/BEF$ & $1=\word{2}$\\
$P_2$ & $BEH/CDG$ & $2=\word{0}$\\
$P_3$ & $AEH/CDF$ & $3=\word{0}$\\\bottomrule
\end{tabular}
\end{center}
Here coordinates are numbered $0,1,2,3$ only in this table.
Their graph is the cycle $P_0P_1P_2P_3P_0$.
The common centers for its four edges are, respectively, $\{\word{1211},\word{1222}\}$, $\{\word{0201},\word{2202}\}$, $\{\word{0000},\word{2100}\}$, and $\{\word{1020},\word{1110}\}$, and the other two pairs have no common center.

A \emph{central hole} of the edge is a hole $q$ that equals $r$ outside $U$ and satisfies $U\subseteq J(q)$.
The contacts of such a $q$ in the four directions of $U$ avoid $r$ outside $U$, so they are eight distinct words of the restriction of the edge, hence all eight words, and each of them agrees with $q$ on $U$ in exactly one coordinate.
At a coordinate $t\in U$ the two contacts in direction $t$ equal $q_t$ and the other six words differ from it, so $q_t$ is the symbol used twice in column $t$ of~\eqref{gen:eight}.
Hence $q$ equals $\word{1200}$ on $U$.
Conversely, $\word{1200}$ extended by $r$ outside $U$ has each of the eight words as a contact, two in each direction of $U$, so it is a central hole.
So the edge has exactly one central hole.
It depends only on $U$ and on the mask outside $U$, and these are the same for both common centers of the edge.

At each object's fixed coordinate, $\word{1200}$ has that object's fixed value.
So the central hole agrees with the outside mask of each of the four objects, and each of their words is a contact of it.
It is therefore one of the three hole representations~\eqref{gen:three-holes} of each of the four objects.
So if two objects are adjacent in $\Omega$ and $q$ is the central hole of their edge, both objects are vertices of $\Gamma_q$, and they are adjacent in $\Gamma_q$ by the definition of $\Gamma_q$.

Conversely, let $q$ be a hole and let $S,T$ be adjacent vertices of $\Gamma_q$, with $U=S\cup T$.
By the proof of~\eqref{gen:edge}, a common center of their objects equals $q$ outside $S\cap T$, and in particular outside $U$.
Since also $U\subseteq J(q)$, the hole $q$ is the central hole of this edge of $\Omega$.
The four three-element subsets of $U$ are vertices of $\Gamma_q$, and their objects are the four objects $P_0,\ldots,P_3$ above, which form a cycle in $\Omega$.
So an edge of $\Gamma_q$ forces all four triples on its four directions to form this square, and all four are nonisolated in $\Gamma_q$.

\subsection{Coherence of the central hole}

\begin{lemma}\label{gen:coherence}
Two edges of $\Omega$ incident to the same object have the same central hole.
\end{lemma}
\begin{proof}
Let $O=(S,m)$ be the common object, and call $S$ the \emph{core}.
Let the two edges join $O$ to objects $O_1$ and $O_2$, and let $q_1$ and $q_2$ be their central holes.
By Section~\ref{gen:square}, $q_1$ and $q_2$ are hole representations of $O$, so both equal $m$ outside $S$.
If the supports of $O_1$ and $O_2$ add the same coordinate $j$ to $S$, then both edges have $U=S\cup\{j\}$ and the mask $m$ outside $U$, so the defining property of their central holes is the same, and $q_1=q_2$.
Suppose they add different coordinates $j$ and $k$, and suppose that $q_1\ne q_2$.
Then $q_1$ and $q_2$ differ on the core.
Normalize $O$ as in~\eqref{gen:six} on the core, and normalize the symbols at $j$ and $k$ so that $m_j=m_k=\word{0}$.
The core projections $h$ and $h'$ of $q_1$ and $q_2$ are then two distinct masks~\eqref{gen:three-holes}.

We use the following symmetries of~\eqref{gen:six}.
Regard a word on the core as a map from its three coordinates to the three symbols.
For permutations $p$ of the symbols and $\pi$ of the coordinates with the same parity, the map $w\mapsto p\circ w\circ\pi^{-1}$ permutes the constant words $A_0,A_1,A_2$ and the even permutation words $B_0,B_1,B_2$, so it preserves~\eqref{gen:six}.
It also permutes the three masks~\eqref{gen:three-holes}, which are the odd permutation words.
These $18$ maps act transitively on the ordered pairs of distinct masks among them.
The stabilizer of the pair $(\word{021},\word{102})$ consists of the identity and the two maps with $p$ a three-cycle and $\pi=p^{-1}$, which permute $A_0,A_1,A_2$ cyclically.
We may therefore assume that $h=\word{021}$ and $h'=\word{102}$, and we may still permute the $A$ labels cyclically without moving $h$ and $h'$.

Call the six words of $O$ the old words.
For an old word $w$ let $f(w)$ and $g(w)$ be its entries at $j$ and $k$.
They are nonzero, since $w$ avoids $m$, and we encode them by signs.
The restriction of the first edge consists of the six old words and two new words $X$ and $Y$.
In it, $O$ is the object with fixed coordinate $j$, whose fixed value $m_j=\word{0}$ is the symbol used twice in column $j$, so $X$ and $Y$ have entry $\word{0}$ at $j$.
They avoid $m$ at $k$, so their entries $x$ and $y$ at $k$ are nonzero.
Likewise, the second edge adds two words $Z$ and $T$ with entry $\word{0}$ at $k$ and nonzero entries $z$ and $t$ at $j$.
The ten words are distinct, since the old words are nonzero at $j$ and $k$, the words $X,Y$ are zero only at $j$, and the words $Z,T$ only at $k$.
All ten avoid $m$ outside the five coordinates $S\cup\{j,k\}$, so their projections to these coordinates are distinct and form a trifferent set.

At $q_1$ the old words are the contact pairs in the directions of the core, and $j\in J(q_1)$.
By the first paragraph of Section~\ref{gen:reps}, the two words of each of these contact pairs have opposite values of $f$.
At $h=\word{021}$ the contact pairs are $(A_0,B_0)$, $(A_2,B_1)$ and $(A_1,B_2)$.
The values $f(A_0),f(A_1),f(A_2)$ are not all equal.
Otherwise $f$ would be constant on $B_0,B_1,B_2$ with the opposite sign, and any two words of the same rainbow triple would agree at $j$.
But the owners of a common center of the first edge are two words of a rainbow triple of $O$, and by Section~\ref{gen:square} they differ at every coordinate of $S\cup\{j\}$.
After exchanging the two nonzero symbols at $j$ if necessary, exactly one of $f(A_0),f(A_1),f(A_2)$ is $-1$, and after a cyclic permutation of the $A$ labels it is $f(A_0)$.

The new words $X$ and $Y$ are the contacts of $q_1$ in direction $j$, so their core prefixes avoid $h$, which leaves eight possible prefixes.
Moreover, let $w,w'$ be old words with $f(w)=f(w')$.
The triple $X,w,w'$ is not rainbow at $j$, where its entries are $\word{0},f(w),f(w')$, nor at $k$, where all three entries are nonzero, nor outside $S\cup\{j,k\}$.
So it must be rainbow on the core, and the same holds for $Y$.
Exactly two of the eight prefixes pass this test for all such pairs $w,w'$, and $X$ and $Y$ have distinct prefixes, since their projections to $S\cup\{j\}$ are distinct and both have entry $\word{0}$ at $j$.
Altogether,
\begin{equation}\label{gen:attachment}
 f(A)=(-,+,+),\quad f(B)=(+,-,-),\quad X=\word{110},\quad Y=\word{202},
\end{equation}
where $X$ and $Y$ now denote the core prefixes.
In the same way, the values $g(A_0),g(A_1),g(A_2)$ are not all equal, and after exchanging the two nonzero symbols at $k$ exactly one of them is $-1$, say $g(A_r)$.
At $h'=\word{102}$ the contact pairs are $(A_1,B_1)$, $(A_0,B_2)$ and $(A_2,B_0)$, so
\[g(B)=(-g(A_2),-g(A_1),-g(A_0)).\]
The same two tests, with $g$ and $h'$ in place of $f$ and $h$, give the core prefixes of $Z$ and $T$ listed in the table below.

Now consider a cross pair $P\in\{X,Y\}$ and $Q\in\{Z,T\}$, and an old word $w$ such that $P,Q,w$ have no rainbow coordinate on the core.
Their rainbow coordinate is then $j$ or $k$.
At $j$ their entries are $\word{0}$, $Q_j$ and $f(w)$, and at $k$ they are $P_k$, $\word{0}$ and $g(w)$.
So
\begin{equation}\label{gen:binary-relations}
 (Q_j,P_k)\ne(f(w),g(w)).
\end{equation}
For each cross pair the table lists two such old words, and~\eqref{gen:binary-relations} for both gives the relations in the last column.
\begin{center}\small
\setlength{\tabcolsep}{3pt}
\begin{tabular}{clll}\toprule
$r$ & $(Z,T)$ & Witness pairs for $XZ,XT,YZ,YT$ & Forced relations\\\midrule
0 & $(\word{220},\word{011})$ & $(A_1,B_1),(A_0,A_1),(A_0,A_2),(B_0,B_2)$
 & $\begin{gathered}z=t=-x,\\z=-y,\ t=y\end{gathered}$\\
1 & $(\word{010},\word{221})$ & $(A_0,A_1),(A_1,B_1),(A_0,B_0),(A_2,B_2)$
 & $\begin{gathered}z=t=x,\\z=y,\ t=-y\end{gathered}$\\
2 & $(\word{211},\word{020})$ & $(A_1,B_1),(A_1,B_1),(A_1,B_2),(A_0,A_2)$
 & $\begin{gathered}z=t=-x,\\z=-y,\ t=y\end{gathered}$\\\bottomrule
\end{tabular}
\end{center}
Every row is contradictory, since its relations give both $x=y$ and $x=-y$.
The prefixes, the witnesses and the contradictions are finite checks, and so are the claims about~\eqref{gen:eight} and its four objects in Section~\ref{gen:square}.
The supplementary script \texttt{anc/general/check\_section3\_tables.py} recomputes all of them.
This exhausts the three possible positions $r$ of the minus sign of $g(A)$ and proves the lemma.
\end{proof}

Consequently an object has a neighbour in $\Gamma_q$ for at most one of its three hole representations $q$.
Indeed, if it has a neighbour in $\Gamma_q$, the corresponding edge of $\Omega$ has central hole $q$ by Section~\ref{gen:square}, and by Lemma~\ref{gen:coherence} all edges of $\Omega$ at the object have the same central hole.
Conversely, if the object has an edge in $\Omega$, it has a neighbour in $\Gamma_q$ for the central hole $q$ of that edge, which is one of its hole representations.
Lemma~\ref{gen:coherence} is the step that relates different holes of the code.
Estimates at each hole separately would not give this conclusion.

\subsection{Counting isolated representations}

A vertex of $\Gamma_q$ is \emph{isolated} if it has no neighbour in $\Gamma_q$.
Let $I(q)$ be the number of isolated vertices of $\Gamma_q$, and put $d=d_q$.
We claim that
\begin{equation}\label{gen:isolated}
 I(q)\le d,\qquad I(q)\le\lfloor d/3\rfloor\quad(d\ge5).
\end{equation}
Only the first bound is used below.
For $d<3$, $d=3$ and $d=4$, the total numbers of triples are $0$, $1$ and $4$.
For $d\ge5$, let $S$ be isolated.
If $a,b\in S$ and $c$ is its third point, nonadjacency of $\{a,b,c\}$ and $\{a,b,x\}$ gives $a_{xa}a_{xb}=-a_{ca}a_{cb}$ for every $x\in J(q)\setminus S$, by~\eqref{gen:edge}.
Hence any two points $x,y\in J(q)\setminus S$ have $a_{xa}a_{xb}=a_{ya}a_{yb}$ for every pair $\{a,b\}\subset S$.
So for every such pair the triples $\{a,b,x\}$ and $\{a,b,y\}$ are adjacent, and their forced square makes all four triples on $\{a,b,x,y\}$ nonisolated.
Another isolated triple cannot meet $S$ in one point $a$: take its other two points as $x,y$ and a second point $b$ of $S$.
It cannot meet $S$ in two points $a,b$ either: take its third point as $x$ and a further point $y\in J(q)\setminus S$, which exists because $d\ge5$.
Thus isolated triples are pairwise disjoint, which proves~\eqref{gen:isolated}.

Let $H_0$ be the number of isolated vertices of $\Omega$, and let $A=\sum_{q\in\F_3^n:\,f(q)=0}d_q$.
Each object is a vertex of $\Gamma_q$ exactly for its three hole representations $q$.
By the end of the previous subsection, an object with an edge in $\Omega$ is isolated in exactly two of these three graphs.
An isolated vertex of $\Omega$ is isolated in all three, since edges of $\Gamma_q$ are edges of $\Omega$.
Therefore, using $I(q)\le d_q$,
\begin{equation}\label{gen:charge}
 \sum_{q\in\F_3^n:\,f(q)=0} I(q)=2(H-H_0)+3H_0=2H+H_0,\qquad 2H\le A.
\end{equation}
There are also exact pair counts:
\begin{align}
 D&=\sum_{\{x,y\}\subset C}2^{n-d(x,y)},\label{gen:D}\\
 A&=\frac12\sum_{\{x,y\}\subset C}(n-d(x,y))2^{n-d(x,y)}.
 \label{gen:A}
\end{align}
For~\eqref{gen:D}, the masks avoided by both words of a pair are double masks with that pair as owners, and every double mask arises in this way from exactly one pair.
For~\eqref{gen:A}, fix a pair $\{x,y\}$ and a coordinate $i$ at which $x$ and $y$ agree.
A mask that equals their common value at $i$ and is avoided by both elsewhere is forced at every coordinate where they differ and has two choices at every other coordinate where they agree.
So there are $2^{n-d(x,y)-1}$ such masks, and each of them is a hole with $i\in J(q)$ whose contact pair in direction $i$ is $\{x,y\}$.
Every pair of a hole $q$ and a direction $i\in J(q)$ arises in this way from exactly one pair, its contact pair in direction $i$, which proves~\eqref{gen:A}.

\subsection{From double boxes to the code size}

For every integer $t\ge0$, $\binom t3\ge3t-8$.
Equality holds at $t=3$ and $t=4$, and from $t$ to $t+1$ the left side grows by $\binom t2$ and the right side by $3$, so the inequality persists for $t\ge4$.
The cases $t\le2$ are immediate.
Summing over the holes, \eqref{gen:objects} gives $3H\ge3A-8Z$, and with~\eqref{gen:charge}
\[
 3A-8Z\le3H\le\tfrac32A,\qquad A\le\tfrac{16}{3}Z.
\]
Now let $M\ge4$, and let $R=d(x,y)$ for a uniformly random pair $\{x,y\}$ of distinct codewords.
By~\eqref{gen:D} and~\eqref{gen:A},
\[
 \frac{2A}{D}=n-\frac{\E[R\,2^{-R}]}{\E[2^{-R}]}.
\]
Since $2^{-R}$ is a decreasing function of $R$, its covariance with $R$ is nonpositive, so $\E[R\,2^{-R}]\le\E[R]\,\E[2^{-R}]$.
In one coordinate, if the three symbols occur $m_0,m_1,m_2$ times, there are $m_0m_1+m_0m_2+m_1m_2\le M^2/3$ pairs of codewords that differ there.
So $\E[R]\le nM^2/3\big/\binom M2=2nM/[3(M-1)]$, and
\[
 A\ge\frac D2\Bigl(n-\frac{2nM}{3(M-1)}\Bigr)=\frac{n(M-3)}{6(M-1)}D\ge\frac n{18}D,
\]
since $(M-3)/(M-1)\ge1/3$ for $M\ge4$.
Combining this with $A\le16Z/3$, without dividing by $D$ or $Z$, we get, for $M\ge4$,
\begin{equation}\label{gen:double-bound}
 D\le\frac{96}{n}Z\le\frac{96}{n}3^n.
\end{equation}

Codes of size $M\le3$, including the empty code, satisfy the theorem directly.
For $M\ge4$ put $\mu=M(2/3)^n$, $\delta=D/3^n$ and $s=2D/(M2^n)$, so that $2\delta=\mu s$.
The final collision estimate is a discrete moment argument, as in \cite[Lemma~3]{DellaFiore2026}.
Let $X,Y$ be independent uniform elements of $C$.
Since $2^{n-d(X,Y)}$ is the number of sets $I\subseteq[n]$ on which $X$ and $Y$ agree, and since the ordered pairs with $X=Y$ contribute $M2^n$ and the others $2D$,
\[
 \E 2^{n-d(X,Y)}=\sum_{I\subseteq[n]}\Pr[X_I=Y_I]
                 =\frac{M2^n+2D}{M^2}=\frac{2^n}{M}(1+s).
\]
The projection $X_I$ takes at most $3^{|I|}$ values, so by the Cauchy--Schwarz inequality $\Pr[X_I=Y_I]\ge3^{-|I|}$.
Also $\Pr[X_I=Y_I]\ge\Pr[X=Y]=1/M$.
Use the first bound when $|I|\le k_0$ and the second for the other sets.
With
\[
 \alpha=\Pr[\operatorname{Bin}(n,1/4)\le k_0],\qquad
 \beta=\Pr[\operatorname{Bin}(n,1/2)\le k_0],
\]
we have $\sum_{|I|\le k_0}3^{-|I|}=(4/3)^n\alpha$, and there are $2^n(1-\beta)$ sets with $|I|>k_0$.
So $(2^n/M)(1+s)\ge(4/3)^n\alpha+2^n(1-\beta)/M$, which gives
\begin{equation}\label{gen:collision}
 s\ge\alpha\mu-\beta,\qquad
 2\delta\ge\alpha\mu^2-\beta\mu.
\end{equation}
Take $k_0=\lfloor n/3\rfloor$.
By Hoeffding's inequality, with deviations at least $n/12$ from the mean $n/4$ and $n/6$ from the mean $n/2$, we have $1-\alpha\le e^{-n/72}$ and $\beta\le e^{-n/18}$.
Solving the quadratic inequality~\eqref{gen:collision} for $\mu$, using $\sqrt{\beta^2+8\alpha\delta}\le\beta+\sqrt{8\alpha\delta}$ and~\eqref{gen:double-bound}, we get
\[
 \mu\le\frac\beta\alpha+\sqrt{\frac{2\delta}{\alpha}}
      \le\frac\beta\alpha+\sqrt{\frac{192}{\alpha n}}
      =O(n^{-1/2}).
\]
This proves Theorem~\ref{thm:half-power}.
Keeping $(M-3)/(M-1)$ in the preceding distance bound gives $D\le32(M-1)3^n/[n(M-3)]$.
If $M\ge n$, the same quadratic estimate yields $\mu\le(8+o(1))/\sqrt n$, uniformly.
If $M<n$, then $\mu\sqrt n\le n^{3/2}(2/3)^n=o(1)$.
Thus the leading constant $8+o(1)$ is valid for arbitrary sequences of codes, although it is not needed for the polynomial improvement.

\phantomsection
\section*{AI use disclosure}
\addcontentsline{toc}{section}{AI use disclosure}
We began by trying to improve the upper bound on the size of a trifferent code by proving that three independent samples from any distribution must fail to separate with at least a certain probability.
The first author pointed out that the search for counterexamples could be turned around: a distribution with a lower failure probability per coordinate than the K\"orner--Marton source would give a better construction by taking independent copies and then deleting words from bad triples.
ChatGPT and Codex helped find and count such examples, search for stronger local rules and mixtures of degrees, and check the resulting estimates in exact arithmetic.
For the overlapping construction, the authors organized the counting argument around a fixed pattern of constraints and clarified why close pairs are removed before the remaining bad triples are counted.
For the upper bound, the authors directed a separate search to improve the polynomial factor in their earlier result.
The models helped develop the key argument that two conflicts involving the same six-word configuration must point to the same empty box.
Separate model runs checked the local cases and searched for counterexamples.
The authors take responsibility for the final proofs and exposition.

\appendix
\section{The finite certificates}
\label{cert:appendix}\label{lower:details}

Section~\ref{lower:distance} reduces the bounds~\eqref{lower:budgets} to finitely many strict inequalities between explicit numbers.
These are the vertex inequalities of Lemma~\ref{lower:cell-lemma} on finitely many cells, and the integer comparison~\eqref{lower:first-moment-check}.
The certificate consists of three kinds of data.
The first is the integer coefficient tables of $H_Q$ and $F_Q$ (Appendix~\ref{lower:local-counts}).
The second is, for pairs and for triples, a subdivision tree of the domain of port profiles, whose leaves are the cells (Appendix~\ref{cert:certificates}).
The third is a rational witness for each retained cell (Appendix~\ref{cert:certificates}).
The verification proceeds in dependency order.
The coefficient tables determine the rates $B_j$ (Appendix~\ref{lower:local-counts}).
The necessary conditions of Appendix~\ref{lower:domain-details} justify the excluded cells.
The subdivision trees prove that the cells cover the domain (Appendix~\ref{cert:certificates}).
Finally, the strict vertex inequalities are proved with rigorous logarithm enclosures (Appendix~\ref{cert:logarithms}).

\subsection{The local data}
\label{lower:local-counts}

For the local rule~\eqref{lower:Q}, write $N_Q=|Q|$.
Then
\[
 N_Q=9\cdot8\frac{14!}{4!\,3!}=43589145600.
\]
The factors $9$ and $8$ choose the labels of multiplicities four and three.
Since $X,Y,Z$ are independent and uniform on $Q$, the local polynomials~\eqref{lower:local-polynomials} are the finite sums
\begin{align*}
 H_Q&=N_Q^{-2}\sum_{x,y\in Q}w_{\mathsf E}^{14-d(x,y)}w_{\mathsf N}^{d(x,y)},\\
 F_Q&=N_Q^{-3}\sum_{\substack{x,y,z\in Q\\(x,y,z)\text{ bad}}}
 \prod_{s\in\Sigma_3} w_s^{N_s(x,y,z)}.
\end{align*}
Permutations of the positions and permutations of the labels map $Q$ onto itself and preserve all states, and they act transitively on $Q$.
So the inner sums over $y$, and over $(y,z)$, do not depend on $x$.
Fix any $x\in Q$.
Then
\begin{equation}\label{cert:count-polynomials}
 H_Q=\frac1{N_Q}\sum_{j=0}^{14}H_jw_{\mathsf E}^jw_{\mathsf N}^{14-j},
 \qquad
 F_Q=\frac1{N_Q^2}\sum_{\substack{\eta\colon\Sigma_3\to\mathbb Z_{\ge0}\\|\eta|=14}}T_\eta\prod_{s\in\Sigma_3}w_s^{\eta_s}.
\end{equation}
Here $H_j$ counts $y\in Q$ with $j$ agreements with $x$, and $T_\eta$ counts ordered pairs $(y,z)\in Q^2$ having no rainbow coordinate with $x$ and equality-state counts $\eta$.

To recover the triple table, append one column at a time and retain only nonrainbow columns.
Number the labels of $\F_9$ as $1,\ldots,9$, and let $\mathbf e_i$ be the $i$-th standard basis vector of $\mathbb Z^9$.
For count vectors $\mathbf a,\mathbf b,\mathbf c\in\mathbb Z_{\ge0}^9$ with a common total, let $T(\mathbf a,\mathbf b,\mathbf c)$ be the generating polynomial of the ordered pairs of second and third rows, with symbol counts $\mathbf b$ and $\mathbf c$, that complete a first row with symbol counts $\mathbf a$ without a rainbow column.
Each column contributes the factor $w_s$ for its state $s$.
Permuting the positions of all three rows together does not change this polynomial, so it depends on the first row only through $\mathbf a$.
Let $i$ be a label with $a_i>0$, and place it in the next column of the first row.
Any such choice gives the same polynomial, and the supplied program chooses one that limits branching.
Choosing the labels $j$ and $k$ of that column in the other two rows gives
\begin{equation}\label{cert:count-recurrence}
 T(\mathbf a,\mathbf b,\mathbf c)
 =\sum_{\substack{j,k\in[9]:\ b_jc_k>0\\|\{i,j,k\}|\le2}}
 w_{\sigma(i,j,k)}
 T(\mathbf a-\mathbf e_i,\mathbf b-\mathbf e_j,\mathbf c-\mathbf e_k),
 \qquad T(0,0,0)=1,
\end{equation}
where $\sigma(i,j,k)\in\Sigma_3$ is the state of the labels $i,j,k$.
Fix $\mathbf a$ to $(4,3,1^7)$ and sum over all distinct assignments of that profile to $\mathbf b$ and $\mathbf c$.
This gives the numerator of $F_Q$.
The order chosen for the first row only permutes coordinates.

The implementation combines count configurations $(\mathbf a,\mathbf b,\mathbf c)$ related by label permutations, and also those related by exchanging the second and third rows, which permutes the pair states accordingly.
For a fixed first-row profile $\mathbf a$, the configurations with a given multiset of triples $(a_i,b_i,c_i)$ form one orbit under the label permutations that preserve $\mathbf a$.
Put
\[
 f_u=|\{i:a_i=u\}|,\qquad
 f_{uvw}=|\{i:(a_i,b_i,c_i)=(u,v,w)\}|.
\]
The label permutations that preserve $\mathbf a$ form a group of order $\prod_{u\ge0}f_u!$, and the stabilizer of $(\mathbf a,\mathbf b,\mathbf c)$ in it has order $\prod_{u,v,w\ge0}f_{uvw}!$.
So the orbit has
\[
 \frac{\prod_{u\ge0}f_u!}{\prod_{u,v,w\ge0}f_{uvw}!}
\]
elements.
There are 87 classes of initial configurations for $Q$, yielding 680 aggregate coefficients.
The supplied recount, run with the option \texttt{--recount}, evaluates every case from the recurrence, checks its multiplicity, and compares every coefficient with the stored table.
Without that option the verifier uses the stored coefficients of each case.
The recount reruns the program that produced the table, so it guards against corrupted data rather than against an error in the program.

The pair counts receive an independent check.
Let $\mathrm{Rearr}(\mathbf a)$ be the set of distinct vectors obtained by permuting the entries of $\mathbf a$, and put
\[
 \Delta_j=(14-j)![u^j]
 \sum_{\mathbf b\in\mathrm{Rearr}(\mathbf a)}\prod_{i=1}^{9}
 \left(\sum_{k=0}^{\min(a_i,b_i)}
       \frac{\binom{a_i}{k}}{(b_i-k)!}u^k\right).
\]
Here $[u^j]$ denotes the coefficient of $u^j$.
This counts a word together with $j$ selected matching positions, so $\Delta_j=\sum_{i=j}^{14}\binom ij H_i$.
Binomial inversion recovers $H_i=\sum_{j=i}^{14}(-1)^{j-i}\binom ji\Delta_j$, since $\sum_{j=i}^{i'}(-1)^{j-i}\binom ji\binom{i'}j=\binom{i'}i(1-1)^{i'-i}$ for $i'\ge i$.
The verifier also checks the collision face $T_{(j,14-j,0,0)}=H_j$ and every permutation of the three $B$-state exponents.
These checks tie the pair and triple tables together.

For completeness, the affine permutations can be written entirely over $\F_3$: use $\F_9=\F_3[\mathrm i]/(\mathrm i^2+\word{1})$, write $\alpha=\alpha'+\alpha''\mathrm i\ne\word{0}$ and $\beta=\beta'+\beta''\mathrm i$, and send
\[
 (a,b)\longmapsto(\alpha'a+\word{2}\alpha''b+\beta',\;\alpha''a+\alpha'b+\beta'').
\]
The determinant is $\alpha'^2+\alpha''^2\ne\word{0}$, giving the $8\cdot9=72$ permutations used in the construction.

\subsection{Safe restrictions on the port profiles}
\label{lower:domain-details}

The weights in~\eqref{lower:tilted-sum} already account for the distance conditions, since they bound the indicator of $\Omega_j$ from above, by~\eqref{lower:weight-domination}.
The following necessary conditions on port counts make the certificate smaller.
They are used only to discard cells that cannot contain a relevant port profile.
They are derived for the certified parameters of Section~\ref{lower:certified}, that is, $m_2=9k$, $m_3=22k$, $D=84k$ and $\delta=3/8$, and their constants depend on these values.

In this subsection, write $\mathsf B_1,\mathsf B_2,\mathsf B_3$ for the pair states $\mathsf B_{12},\mathsf B_{13},\mathsf B_{23}$, and write the port profile of a triple as $h=(h_{\mathsf A},h_1,h_2,h_3)$, where $h_{\mathsf A}$ is the all-equal fraction and $h_i=h_{\mathsf B_i}$ for $i=1,2,3$.
Every far bad triple satisfies
\begin{equation}\label{lower:domain}
 h_i+h_{i'}\ge\frac{69}{224}\ (i\ne i'),\qquad
 h_{\mathsf A}\le\frac{31}{64},\qquad
 2h_{\mathsf A}-3h_i\le\frac{31}{56}\ (i=1,2,3).
\end{equation}
Here is a direct proof.
A set of $u$ variables occupies at least $3u-m_2$ ports, since at most $m_2$ of them can have degree two.
Each pair differs in at least $\delta m$ variables, so its port count is at least $3\delta m-m_2=69D/224$.
Every non-all-equal variable contributes two to the sum of the three pair distances.
There are therefore at least $3\delta m/2$ such variables, occupying at least $9\delta m/2-m_2=33D/64$ ports.
This proves the first two conditions.

For the last, let $u_{\mathsf A},u_i$ be the numbers of variables in states $\mathsf A,\mathsf B_i$, and let $u_{d,s}$ be the number of degree-$d$ variables in state $s$.
Adding the two distance inequalities that involve state $\mathsf B_i$ gives $u_{\mathsf A}-u_i\le m(1-2\delta)$.
Consequently
\[
 D(2h_{\mathsf A}-3h_i)
 =6(u_{\mathsf A}-u_i)-2u_{2,\mathsf A}-3u_{3,i}
 \le6m(1-2\delta)=31D/56.
\]
For a close pair, at most three ports occur per differing variable, so its unequal port fraction satisfies
\begin{equation}\label{lower:pair-domain}
 h_{\mathsf N}\le 3\delta m/D=93/224.
\end{equation}

\subsection{The two certificates for the code construction}
\label{cert:certificates}

This subsection describes how Proposition~\ref{cert:lower-norms} is verified.
The parameters are those of Section~\ref{lower:distance}, all rational: $g=1162579/1000000$, $\alpha_2=3/14$, $\alpha_3=11/14$, $e=31/84$ and $\delta=3/8$, so that $\kappa_2=-4e\log g$ and $\kappa_3=-8e\log g$.
The strict inequalities $B_j<\kappa_j$ are checked at the vertices of each retained cell, with the witness of that cell.
The proof of Proposition~\ref{lower:budgets-prop} combines these inequalities with the coverage, exclusion, tilt and symmetry checks described below.
For pairs, the cover starts from the closed interval with endpoints $(1,0)$ and $(0,1)$.
For triples, $F_Q$ is symmetric in the three pair states (proof of Proposition~\ref{lower:budgets-prop}), and the verifier checks this symmetry on the stored coefficients as a consistency test.
By the symmetry argument in the proof of Proposition~\ref{lower:budgets-prop}, it therefore suffices to cover the ordered region $h_{\mathsf B_{12}}\ge h_{\mathsf B_{13}}\ge h_{\mathsf B_{23}}$, whose closed simplex has vertices
\[
 (1,0,0,0),\quad(0,1,0,0),\quad
 (0,\tfrac12,\tfrac12,0),\quad
 (0,\tfrac13,\tfrac13,\tfrac13).
\]
Repeated midpoint bisections of simplex edges produce rational vertices.
The replay checks that the supplied leaves form a complete binary subdivision: no leaf is an ancestor of another, and every internal node has both children from the same edge bisection.
This proves coverage, including every boundary.
A cell is excluded only if there is one condition among~\eqref{lower:domain}, for triples, or~\eqref{lower:pair-domain}, for pairs, that every vertex of the cell violates strictly.
The verifier also checks that every witness has positive entries and that its tilt satisfies $t\le1$ for pairs and $t_{12},t_{13},t_{23}\ge1$ for triples.
The distance weights remain in $B_j$ even on retained cells.
The cuts merely avoid checking regions that cannot occur.

The pair partition has ten leaves: one exclusion and nine retained cells, requiring eighteen vertex comparisons.
The triple partition has 11,713 leaves: 300 exclusions and 11,413 retained cells, requiring 45,652 vertex comparisons.
These are exhaustive finite checks, not samples from the domain.
The verifier also checks the uniform margins
\[
 B_2+4e\log g<-10^{-5},\qquad
 B_3+8e\log g<-10^{-8}
\]
on every retained cell.
The largest values of the left sides over all vertices of retained cells are about $-1.73\cdot10^{-5}$ and $-1.28\cdot10^{-8}$.
It also checks the integer comparison~\eqref{lower:first-moment-check}.

\subsection{How the logarithms are verified exactly}
\label{cert:logarithms}

At a rational vertex with rational witnesses, each comparison in Proposition~\ref{cert:lower-norms} is a sum of rational multiples of logarithms of positive rational numbers.
Polynomial evaluations are exact.
For a positive rational $\xi$, write
\[
 \xi=2^py,\qquad p\in\mathbb Z,\quad 1\le y<2,
 \qquad\text{so that}\qquad
 \log\xi=p\log2+\log y,
\]
and set
\[
 \chi=\frac{y-1}{y+1}\in\Bigl[0,\frac13\Bigr),\qquad\text{so that}\qquad y=\frac{1+\chi}{1-\chi}.
\]
For every $\chi\in[0,1)$, with $y=(1+\chi)/(1-\chi)$, and every number $q\ge1$ of terms,
\begin{equation}\label{cert:log-enclosure}
 L_q(\chi)\le\log y\le L_q(\chi)+\frac{2\chi^{2q+1}}{(2q+1)(1-\chi^2)},
 \qquad\text{where}\qquad
 L_q(\chi)=2\sum_{i=0}^{q-1}\frac{\chi^{2i+1}}{2i+1}.
\end{equation}
Indeed, the geometric identity
\[
 \frac1{1-x^2}=\sum_{i=0}^{q-1}x^{2i}+\frac{x^{2q}}{1-x^2}
\]
gives
\[
 \log y=2\int_0^\chi\frac{dx}{1-x^2}
 =L_q(\chi)+2\int_0^\chi\frac{x^{2q}}{1-x^2}\,dx,
\]
and, since $1-x^2\ge1-\chi^2$ for $0\le x\le\chi$,
\[
 0\le2\int_0^\chi\frac{x^{2q}}{1-x^2}\,dx\le\frac2{1-\chi^2}\int_0^\chi x^{2q}\,dx=\frac{2\chi^{2q+1}}{(2q+1)(1-\chi^2)}.
\]
Since $2=(1+\frac13)/(1-\frac13)$, the bounds~\eqref{cert:log-enclosure} with $\chi=1/3$ enclose $\log2$.
Combining the two enclosures gives an enclosure of $\log\xi$, and endpoints are reversed when multiplying by a negative coefficient.
The remainder term in~\eqref{cert:log-enclosure} increases with $\chi$, so for $0\le\chi\le1/3$ it is at most its value at $\chi=1/3$, the uniform tail bound
\[
 \frac{2\cdot3^{-(2q+1)}}{(2q+1)\cdot\frac89}=\frac{9}{4(2q+1)3^{2q+1}}.
\]
The replay uses $q=80$ terms and integer fixed-point intervals of scale $2^{160}$, rounding every operation outwards.
A negative upper endpoint for $B_2+4e\log g$ or $B_3+8e\log g$ proves the required strict comparison.
Decimal strings in the witness files are interpreted as exact rationals.
No floating-point value enters a pass/fail decision.

\bibliographystyle{plain}
\bibliography{references}

\begin{thebibliography}{1}

\bibitem{BhandariKhetan2025}
Siddharth Bhandari and Abhishek Khetan.
\newblock Improved upper bound for the size of a trifferent code.
\newblock {\em Combinatorica}, 45:2, 2025.
\newblock Preprint \href{https://arxiv.org/abs/2402.02390}{arXiv:2402.02390}.

\bibitem{DellaFiore2026}
Stefano Della~Fiore.
\newblock A note on vector trifferent codes over the sphere.
\newblock \href{https://arxiv.org/abs/2608.04256v1}{arXiv:2608.04256v1}, 2026.

\bibitem{HanSunGe2026}
Xuejiao Han, Yubo Sun, and Gennian Ge.
\newblock Polynomially improved lower bounds for trifferent codes via locally
  sparse $3$-uniform hypergraphs.
\newblock \href{https://arxiv.org/abs/2607.26376v2}{arXiv:2607.26376v2}, 2026.

\bibitem{KornerMarton1988}
J{\'a}nos K{\"o}rner and Katalin Marton.
\newblock New bounds for perfect hashing via information theory.
\newblock {\em European Journal of Combinatorics}, 9(6):523--530, 1988.

\bibitem{YuanZhu2026}
Chen Yuan and Ruiqi Zhu.
\newblock Beating the probabilistic bound for hashing and improving lower
  bounds for trifferent codes.
\newblock \href{https://arxiv.org/abs/2611.11262}{arXiv:2611.11262}, 2026.

\end{thebibliography}
\end{document}